\documentclass[conference]{IEEEtran}

\usepackage{url}

\usepackage[]{fancyhdr} %
\newcommand{\changefont}{\fontsize{9}{9}\selectfont}
\usepackage{mathptmx}       
\usepackage{helvet}     
\usepackage{courier}        
\usepackage{type1cm}        
\usepackage{makeidx}         
\usepackage{graphicx}        
         
\usepackage{multicol}        
\usepackage[bottom]{footmisc}
\usepackage{ifthen}
\usepackage{subfigure}
\usepackage[usenames,dvipsnames]{xcolor}
\usepackage{pdflscape} 

\makeindex             

\usepackage{graphics} 
\usepackage{graphicx} 
\usepackage{epsfig} 
\usepackage{mathptmx} 
\usepackage{times} 
\usepackage{mathtools}  
\usepackage{amssymb}  
\usepackage{amsthm}

\usepackage{amsfonts}

\usepackage{subfig}
\usepackage{verbatim}
\usepackage{comment}
\usepackage[ruled,vlined]{algorithm2e}	
\usepackage{algorithmic}
\usepackage{multirow}

\def\ba{\begin{array}}
\def\ea{\end{array}}
\newcommand{\beq}{\begin{equation}}
\newcommand{\eeq}{\end{equation}}
\newcommand{\bq}{\begin{eqnarray}}
\newcommand{\eq}{\end{eqnarray}}
\newcommand{\bqn}{\begin{eqnarray*}}
\newcommand{\eqn}{\end{eqnarray*}}
\newcommand{\bee}{\begin{enumerate}}
\newcommand{\eee}{\end{enumerate}}
\newcommand{\bi}{\begin{itemize}}
\newcommand{\ei}{\end{itemize}}
\newcommand{\ii}{\mathbf{i}}

\newcommand{\diag}{\mathrm{diag}}

\newboolean{showcomments} 
 \setboolean{showcomments}{false}
\newcommand{\slow}[1]{\ifthenelse{\boolean{showcomments}}
{ \textcolor{red}{(SL: #1)}}{}}
\newcommand{\czhao}[1]{\ifthenelse{\boolean{showcomments}}
{ \textcolor{ForestGreen}{(CZ: #1)}}{}}

\newboolean{showToDo}
 \setboolean{showToDo}{false}
\newcommand{\ToDo}[1]{\ifthenelse{\boolean{showToDo}}
{ \textcolor{red}{#1}}{}}

\newboolean{showSolution}
 \setboolean{showSolution}{true}
\newcommand{\incSolution}[1]{\ifthenelse{\boolean{showSolution}}
{#1} {} }

\newboolean{showArchive}
 \setboolean{showArchive}{false}
\newcommand{\slowArchive}[1]{\ifthenelse{\boolean{showArchive}}
{#1} {} }

\theoremstyle{plain}
\theoremstyle{definition}
\newtheorem{theorem}{Theorem}
\newtheorem{lemma}[theorem]{Lemma}

\theoremstyle{definition}

\newtheorem{remark}{Remark}

\newcounter{AssumptionCounter}
\newcommand{\AssumptionInc}[1]{\refstepcounter{AssumptionCounter}\label{#1}}

\begin{document}

%

%
 \title{A Three-phase Power Flow Model and Balanced Network Analysis}

\author{\IEEEauthorblockN{Steven H. Low \\ CMS, EE, Caltech}}


%





\maketitle
\thispagestyle{fancy}
\pagestyle{fancy}


\begin{abstract}
First we present an approach to formulate unbalanced three-phase
power flow problems for general networks that explicitly separates 
device models and network models. 
A device model consists of 
(i) an internal model and (ii) a conversion rule.  The conversion rule
relates the internal variables (voltage, current, and power) of a device 
to its terminal variables through a conversion matrix $\Gamma$ and 
these terminal variables are related by network equations.
Second we apply this approach to balanced three-phase networks to formalize 
per-phase analysis and prove its validity for general networks using the 
spectral property of the conversion matrix $\Gamma$.
\end{abstract}

\begin{IEEEkeywords}
Unbalanced three-phase power flow models, balanced networks, per-phase analysis.
\end{IEEEkeywords}


%
\IEEEpeerreviewmaketitle


\section{Introduction}

\noindent\textbf{Motivation.}
Unbalanced three-phase load flow problems are becoming increasingly important as we
decarbonize the grid.  Such problems are more difficult for several reasons; see, e.g., 
\cite[Chapter 11]{GomesExposito2018} for transmission systems and \cite{Kersting2002}
for distribution systems.
First, a network model is more complicated because 
three-phase lines couple currents and voltages in different phases when lines are not 
transposed or loads are unbalanced, e.g., as in most distribution systems.
If the network is symmetric, then a similarity transformation due to Fortescue \cite{Fortescue1918} 
from the phase coordinate to a sequence coordinate produces network models that are
decoupled in the sequence coordinate and can therefore be analyzed in a way similar to
a single-phase network; see, e.g., \cite{Clarke1943}.   Without symmetry, this transformation
however offers no simplification.
Second, the voltages and currents across the single-phase 
devices internal to $\Delta$ configuration are observed externally only through a
linear map $\Gamma$  that is not invertible.  While we are typically interested in solving for or
optimizing the \emph{internal} currents or power flows across the single-phase devices, e.g.,
controlling the charging currents of electric vehicle chargers in $\Delta$ configuration,
a network model, such as $I = YV$ or 
$s_j = \sum_k y_{jk}^s \left(|V_j|^2 - V_jV_k^{\sf H} \right)$, relates only the \emph{terminal}
voltages and currents observable externally of three-phase devices.  The interplay 
between internal and external variables of a three-phase device sometimes seems confusing.
Third, load flow formulations sometimes implicitly assume that the neutrals of all $Y$-configured
devices are at zero potential and the zero-sequence components of the terminal voltages
of all $\Delta$-configured devices are zero.  This limits their applicability as, e.g., they exclude
the case where some $Y$-configured loads are ungrounded or grounded with nonzero 
earthing imepdances.
Solutions to the last two difficulties both lie in a careful accounting of the 
conversion between internal and external variables of $\Delta$-configured devices, 
using the conversion matrix $\Gamma$.

\noindent\textbf{Summary.}
In this paper we present such a modeling approach that separates three-phase device
models and network networks.  In this approach,  a device model consists of two components:
(i) an internal model and (ii) a conversion rule.  The internal model 
describes how each of the single-phase device behaves regardless of their configuration.
The conversion rule, on the other hand, depends only on their configuration regardless of the
type of devices.  It maps internal voltages, currents,
and powers across these single-phase devices to terminal voltages, currents, and powers
observable externally.  
Since the network model relates only the terminal variables regardless of the type of
devices or their configurations, the explicit separation of device and network models 
allows mix and match of equivalent models, enhancing modeling flexibility.  
We present our model in Section \ref{sec:3pModel} and use it to formulate
an unbalanced three-phase analysis problem for general networks in Section \ref{sec:3pAnalysis}.

In Section \ref{part:networks; ch:mun.bim; sec:BalancedNk} we illustrate our model by 
showing formally that the analysis problem can be solved using per-phase analysis when
the network is balanced.   It is well known
that a balanced three-phase device in $\Delta$ configuration has an $Y$ equivalent 
that has the same external behavior.  The standard way to justify per-phase analysis is 
by analyzing specific three-phase circuits, often simple circuits; e.g., 
\cite{ChuaDesoerKuh, Bergen2000, GloverSarmaOverbye2008}, by first converting all
 $\Delta$-configured devices into their $Y$ equivalents, and then showing that all neutrals
in the equivalent circuit are at the same potential and that all phases are decoupled.
This implies that the original three-phase circuit can be solved by analyzing a simpler 
per-phase circuit.  This process has two limitations.  

First, it is not clear how to extend circuit analysis methods, e.g., loop analysis or 
mesh analysis \cite[Chapter 12]{ChuaDesoerKuh}, from specific (and simple) circuits to 
an arbitrary balanced network and prove that all neutrals are at the same potential.
This is simple to show, however, in our model
by expressing the network equation $I = YV$ in terms of the Kronecker product and using the
spectral property of $\Gamma$ (see Theorems \ref{ch:mun.1; sec:bim; subsec:BalancedNk; thm:1phiNk.C1C2}
and \ref{ch:mun.1; sec:bim; subsec:BalancedNk; thm:1phiNk}).
The intuition is as follows.
In a balanced three-phase network, positive-sequence voltages and 
currents are in span$(\alpha_+)$ where $\alpha_+$ 
is an eigenvector of $\Gamma$ and $\Gamma^{\sf T}$.  This means that the 
transformation of balanced voltages and currents under $\left( \Gamma, \Gamma^{\sf T}\right)$ reduces 
to a scaling of these variables by their eigenvalues $1-\alpha$ and $1-\alpha^2$ respectively.  
The voltage and current at every point in a network can be written
as linear combinations of transformed source voltages and source currents,
transformed by $\left( \Gamma, \Gamma^{\sf T} \right)$ and line admittance matrices.
Therefore if the source voltages and source currents are in span$(\alpha_+)$
and if lines are identical and phase-decoupled, then
the transformed voltages and currents remain in span$(\alpha_+)$ and hence 
are balanced positive-sequence sets.  This is explained in 
Sections \ref{part:networks; ch:mun.bim; sec:BalancedNk; subsec:1pNetwork}, 
\ref{part:networks; ch:mun.bim; sec:BalancedNk; subsec:1pAnalysis} and 
\ref{part:networks; ch:mun.bim; sec:BalancedNk; subsec:equivalence}.

Second, all neutrals are at the same potential only if the neutral voltages of all $Y$-configured
devices (voltage sources, current sources and impedances) are assumed zero and the 
zero-sequence voltages of all $\Delta$-configured voltage sources are assumed zero.  Roughly,
this requires that all neutrals are grounded directly (i.e., with zero earthing impedances).  The
standard analysis of specific circuits often makes this assumption sometimes implicitly.  
Without this assumption, the neutral voltages on the circuit are generally different.  Yet, 
per-phase analysis can be extended to the general case without this assumption 
as long as the network is balanced, except that per-phase analysis is needed not only on 
a per-phase positive-sequence network, but also on a per-phase zero-sequence network.  
This is explained in Section \ref{part:networks; ch:mun.bim; sec:BalancedNk; subsec:extension}. 

We close our summary with two remarks on per-phase analysis.  
First, if the network is unbalanced but symmetric,
i.e., impedances are balanced and lines are (phase coupled but) symmetric, then Fortescue's 
similarity transformation \cite{Fortescue1918} from the phase coordinate to the sequence 
coordinate leads to decoupled device models and network models.   
The network equation is therefore decoupled in the sequence coordinate and can be interpreted 
as defining three separate sequence networks, to which the per-phase analysis of
Section \ref{part:networks; ch:mun.bim; sec:BalancedNk} can be applied.
Second, with today's abundant computing power the smaller problem size may not be an 
important advantage of per-phase analysis.  Rather, per-phase analysis illustrates the
application of our modeling approach to three-phase power flow.  It also clarifies 
the simple structure underlying a balanced network and enhances our conceptual
understanding of three-phase networks in general, balanced or unbalanced.

\noindent\textbf{Literature.}
Single-phase models are a good approximation for many transmission network
applications where lines are symmetric and loads are balanced so that the zero and
negative-sequence components are negligible compared with the positive-sequence component.
They may not be negligible when
lines are not transposed or equally spaced, e.g., as in distribution systems, and when loads
are unbalanced or nonlinear, e.g., AC furnaces, high-speed trains, power electronics, or single
or two-phase laterals in distribution networks.  This can cause power quality issues such as 
voltage imbalances and harmonics.  Furthermore single-phase analysis can produce incorrect 
power flow solutions.   

There is a large literature on three-phase power flow analysis and we only make a few
brief remarks.
Three-phase load flow solvers have been developed since at least the 1960s, e.g., see 
\cite{ElAbiadTarsi1967} for solution in the sequence coordinate and \cite{Berg1967, Laughton1968}
in the phase coordinate.   A three-phase network is equivalent to a single-phase circuit
where each node in the equivalent circuit is indexed by a (bus, phase) pair \cite{Laughton1968}.   
Single-phase power flow algorithms such as Newton Raphson \cite{Birt1976} or Fast Decoupled 
methods \cite{ArrillagaArnold1978}
can be directly applied to the equivalent circuit.  The main difference with a single-phase network
is the circuit models of three-phase devices in the equivalent circuit, such as models for
three-phase lines \cite{Chen1991, Kersting2002}, transformers and 
co-generators \cite{Laughton1968, Chen1991Transformer}, 
constant-power devices \cite[Chapter 11]{GomesExposito2018}, as well as 
voltage regulators, and loads \cite{Kersting2002}, etc. 
A state-of-the-art algorithm
in \cite[Chapter 11]{GomesExposito2018} expresses currents in terms of voltages for both $PQ$ 
and $PV$ buses, applies the Newton-Raphson algorithm to the resulting nonlinear current
balance equation $I = YV$ in the sequence domain.  It allows both grounded and ungrounded
loads in $Y$ and $\Delta$ configurations.  For transmission networks, computing in the sequence
domain has the advantage that, when most lines in the network are symmetric and thus have 
decoupled representation in the sequence coordinate, the Jacobian matrix is sparse.  
Sometimes an approximate solution is computed by ignoring the coupling across zero, positive, 
and negative-sequence variables and solving the three sequence networks separately as 
single-phase networks, e.g., \cite{Zhang1996}.
Distribution networks usually does not enjoy such simplification and hence computation is
usually done in the phase coordinate.

While the papers above study general networks that may contain cycles, another set of power 
flow methods are tailored for three-phase radial networks 
 \cite{Berg1967, Kersting1976, Kersting2002, Cheng1995, ZimmermanChiang1995, MiuChiang2000}. 
In particular, the tree topology leads to a spatially recursive structure that enables iterative algorithms
called backward forward sweep (BFS), apparently first developed in \cite{Berg1967}.
Different BFS algorithms are developed in \cite{Kersting1976}\cite[Chapter 10.1.3]{Kersting2002}
\cite{Cheng1995} for three-phase networks (\cite{Cheng1995} generalizing the BFS algorithm of
\cite{Shirmohammadi1988} from single-phase to three-phase networks).  
%
For single-phase radial networks, a solution method based on the DistFlow model is
	developed in \cite{Baran1989b} that uses one-time forward sweep (to express
	all variables in terms of the voltages at the feeder head and all branch points)
	followed by a Newton-Raphson algorithm to solve for these voltages.
	By exploiting the approximate sparsity of the Jacobian 
	matrix in \cite{Baran1989b}, approximate fast decoupled methods are developed 
	and their convergence properties analyzed in \cite{Chiang1991}.
These methods are extended to three-phase radial networks in 
\cite{ZimmermanChiang1995}.  The existence and uniqueness of power flow solutions
of three-phase DistFlow model is analyzed in \cite{MiuChiang2000}.
The advantage of BFS is that it does not need to compute Jacobian nor solve a linear system 
to compute iteration updates.  Newton-Raphson, on the other hand, tends to converge in a smaller
number of iterations.

\slow{
Check out other references for Multiphase radial networks:
\bi
\item Recent papers by Changhong-NREL, Ahmed Zamzam-NREL; 	Mario and Le Boudec's group;
	Federik Geth of CSIRO 3-phase with shunt admittance matrices.
\item Qiuyu's thesis (Chapter 3.2) has a nice description and illustration of notations for unbalanced BFM.
\item Elizondo, Tuffner and Schneider, Three-phase unbalanced transient dynamics and power flow
for modeling distribution systems with synchronous machines, TPS, Jan 2016; and refs [10-14] therein.
\item Emiliano 2013; Gan 2014; Changhong 2017, more...
\ei
}

\noindent\textbf{Notation.}
Let $\mathbb{C}$ denote the set of complex numbers.
For $a\in \mathbb C$, Re$\,\, a$ and Im$\,\, a$ denote its real and imaginary 
parts respectively, and $\bar a$ or $a^{\sf H}$ denotes its complex conjugate.
We use $\ii$ to denote $\sqrt{-1}$.
A vector $x\in \mathbb C^n$ is a column vector and is denoted in one of two ways:
\begin{align*}
x & \ \ = \ \ \begin{bmatrix} x_1 \\ \vdots \\ x_n \end{bmatrix}, 	&
x & \ \ = \ \ (x_1, \dots, x_n) 
\end{align*}
Its componentwise complex conjugate is denoted by $\overline x$.
For any matrix $A$, $A^{\sf T}$, $A^{\sf H}$, $A^\dag$ denote
its transpose, Hermitian transpose, and pseudo-inverse respectively.
If $x$ is a matrix then $\diag\left(x \right)$ is the vector whose components are the diagonal entries
of $x$, whereas if $x$ is a vector then $\diag(x)$ is a diagonal matrix with $x_i$ as its diagonal entries.  
Finally $\textbf 1 \in\mathbb C^3$ is the column vector of size 3 whose
entries are all 1s and $\mathbb I$ is the identity matrix of size 3.

\section{Three-phase network model}
\label{sec:3pModel}

\begin{figure}[htbp]
\centering
   \includegraphics[width=0.4\textwidth] {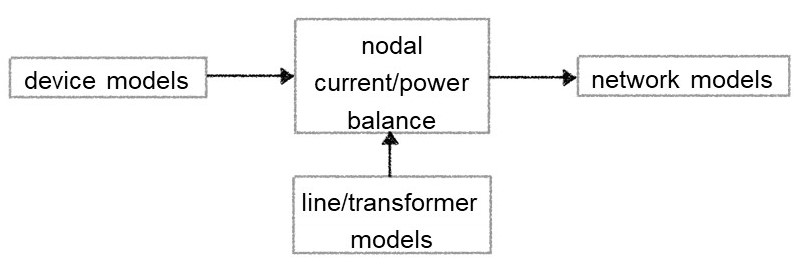}
\caption{Overall three-phase network model.}
\label{fig-IREP-Overall3pModelingSteps}
\end{figure}

A three-phase network connects generators to loads, each of which is modeled by
a single-terminal device.  
Each terminal has three wires (or ports or conductors) indexed
by its phases $a, b, c$, and possibly a neutral wire indexed by $n$.\footnote{For notational simplicity, we assume all devices and lines have three phases. Generalization to the case where 
some devices or lines have only one or two phases is straightforward. }  
Internally, the device can be in 
$Y$ or $\Delta$ 
configuration, and the $Y$ configuration may have a neutral wire that may be grounded.  
A three-phase line has two terminals, each terminal with three or four wires, and it connects two single-terminal devices, one at each end of the line.  Its neutral wire
may be grounded at regular spacing along the line.
The overall network model consists of three components (see Figure 
\ref{fig-IREP-Overall3pModelingSteps}):
\bee
\item \emph{Device model.} The internal behavior of a
    single-terminal device is defined by the relationship between the
	voltages $V^Y/V^\Delta$, currents $I^Y/I^\Delta$, and powers $s^Y/s^\Delta$ across
	each of the single-phase devices that make up the three-phase device.   This relationship
	is independent of whether the device is in $Y$ or $\Delta$ configuration. The configuration
	defines a conversion rule that maps internal
	variables $\left(V^{Y/\Delta}, I^{Y/\Delta}, s^{Y/\Delta}\right) \in \mathbb C^{9}$ 
	to terminal voltages, currents, and powers $(V, I, s) \in \mathbb C^{9}$,  
	regardless of the type of the devices.
	The internal behavior and the conversion rule jointly determine the external behavior of the 
	three-phase device, i.e., the relationship between the terminal variables $(V, I, s)$ that can 
	be observed externally.  
    
\item \emph{Line model.} It relates the terminal voltages, 
    line currents, and line power flows
    $(V_j, I_{jk}, S_{jk}) \in \mathbb C^{9}$ and $(V_k, I_{kj}, S_{kj}) \in \mathbb C^{9}$, at each end of the line $(j,k)$.
    
\item \emph{Network model.} It relates the terminal variables of all devices on the network.
	This is defined by current or power balance at every bus in the network.
\eee

In this section we present circuit models of four types of single-terminal devices, a voltage
source, a current source, a power source, and an impedance.  Then we present a 
circuit model of a three-phase line.  Finally we compose an overall model of
a network of such devices.

\subsection{Conversion matrices $\Gamma, \Gamma^{\sf T}$}
\label{sec:3pBFM; subsec:ConversionMatrices}

We start by defining \emph{conversion matrices} $\Gamma, \Gamma^{\sf T}$ that
maps between internal and external variables in a $\Delta$ configuration:
\begin{subequations}
\bq
\Gamma &  :=  & \begin{bmatrix}  1 & -1 & 0 \\  0 & 1 & -1  \\  -1 & 0 & 1  \end{bmatrix}
\\
\Gamma^{\sf T}&  :=  & \begin{bmatrix}  1 & 0 & -1 \\  -1 & 1 & 0  \\  0 & -1 & 1  \end{bmatrix}
\eq
As we will see,
\label{sec:3pBFM; subsec:ConversionMatrices; eq:defGamma}
\end{subequations}
the spectral properties of $\Gamma, \Gamma^{\sf T}$ underlie much of the behavior
of three-phase systems, balanced or unbalanced.  Here we recall some basic
facts on $\Gamma, \Gamma^{\sf T}$ that are useful in the rest of the paper.

It can be shown that $\Gamma$ and $\Gamma^{\sf T}$ are normal matrices and 
their spectral decompositions are
\begin{align*}
\Gamma & \ \ = \ \ F \Lambda \overline F, 	\qquad\qquad  	\Gamma^{\sf T} \ \ = \ \ \overline F \Lambda F
\end{align*}
where $\Lambda$ is a diagonal matrix and $F$ is a unitary matrix defined as:
\begin{align}
\Lambda & \ \ := \ \ \begin{bmatrix} 0 &  &  \\  & 1-\alpha &  \\  &  & 1-\alpha^2 \end{bmatrix}, 	&
F & \ \ := \ \ \frac{1}{\sqrt 3} \begin{bmatrix} \textbf 1 & \alpha_+ & \alpha_- \end{bmatrix}
\label{part:NetworkModels; ch:mun.1; sec:ComponentModels; subsec:GenLoadD; eq:defF}
\end{align}
with $\alpha := e^{-\ii 2\pi/3}$.  The positive-sequence and
negative-sequence vectors $\left( \alpha_+, \alpha_- \right)$  are
\bqn
\alpha_+ & := & \begin{bmatrix} 1 \\ \alpha \\ \alpha^2 \end{bmatrix}, 	\qquad\qquad
\alpha_-  \ \ := \ \ \begin{bmatrix} 1 \\ \alpha ^2\\ \alpha \end{bmatrix} 
\eqn
The eigenvectors of $\Gamma$ are $\textbf 1, \alpha_+, \alpha_-$ and they are
orthogonal.
Here $\overline F$ is the complex conjugate of $F$ componentwise.
Since $F$ is symmetric, the pseudo inverses of $\Gamma, \Gamma^{\sf T}$ are  
\begin{align*}
\Gamma^\dag & \ \ = \ \ F \Lambda^\dag \overline F, 	\qquad\qquad
\Gamma^{\sf T\dag} \ \ = \ \ \overline F \Lambda^\dag F
\end{align*}
where $\Lambda^\dag := \diag\left(0, (1-\alpha)^{-1}, (1-\alpha^2)^{-1} \right)$.
This yields the following properties.
\begin{lemma}[Pseudo inverses of $\Gamma, \Gamma^{\sf T}$]
\label{part:NetworkModels; ch:mun.1; sec:ComponentModels; subsec:GenLoadD; thm:InverseG}
\bee
\item The null spaces of $\Gamma$ and $\Gamma^{\sf T}$ are both span$(\textbf 1)$.
\item Their pseudo-inverses are
	\bqn
	\Gamma^\dag \ \ = \ \ \frac{1}{3} \, \Gamma^{\sf T},	\qquad\quad
	\Gamma^{\sf T \dag} \ \ = \ \ \frac{1}{3}\,  \Gamma 
	\eqn
\item Consider $\Gamma x = b$ where $b, x\in\mathbb C^3$.  
	Solutions $x$ exist if and only if $\textbf 1^{\sf T} b = 0$, in which case the solutions $x$ are given by 
	\bqn
	x & = & \frac{1}{3}\, \Gamma^{\sf T} b \ + \ \gamma\textbf 1, 	\qquad\qquad \gamma\in\mathbb C
	\eqn
\item Consider $\Gamma^{\sf T} x = b$ where $b, x\in\mathbb C^3$.   
	Solutions $x$ exist if and only if $\textbf 1^{\sf T} b = 0$, in which case the solutions $x$ are given by 
	\bqn
	x & = & \frac{1}{3}\, \Gamma b \ + \ \gamma\textbf 1, 	\qquad\qquad \gamma\in\mathbb C
	\eqn

\item $\Gamma\Gamma^\dag = \Gamma^\dag \Gamma = 
	\frac{1}{3}\, \Gamma\Gamma^{\sf T} \ = \ \frac{1}{3}\,\Gamma^{\sf T} \Gamma \ \ = \ \ \mathbb I \ - \ \frac{1}{3}\, \textbf 1 \textbf 1^{\sf T}$
	where $\mathbb I$ is the identity matrix of size 3.
\eee
\end{lemma}
In this paper, we will use Lemma 
\ref{part:NetworkModels; ch:mun.1; sec:ComponentModels; subsec:GenLoadD; thm:InverseG} repeatedly, sometimes without explicit reference.

\subsection{Devices: internal models and conversion rules}
\label{sec:3pBFM; subsec:InternalModels}

\begin{figure}[htbp]
\centering
\subfigure [$Y$ configuration]{
   \includegraphics[width=0.35\textwidth] {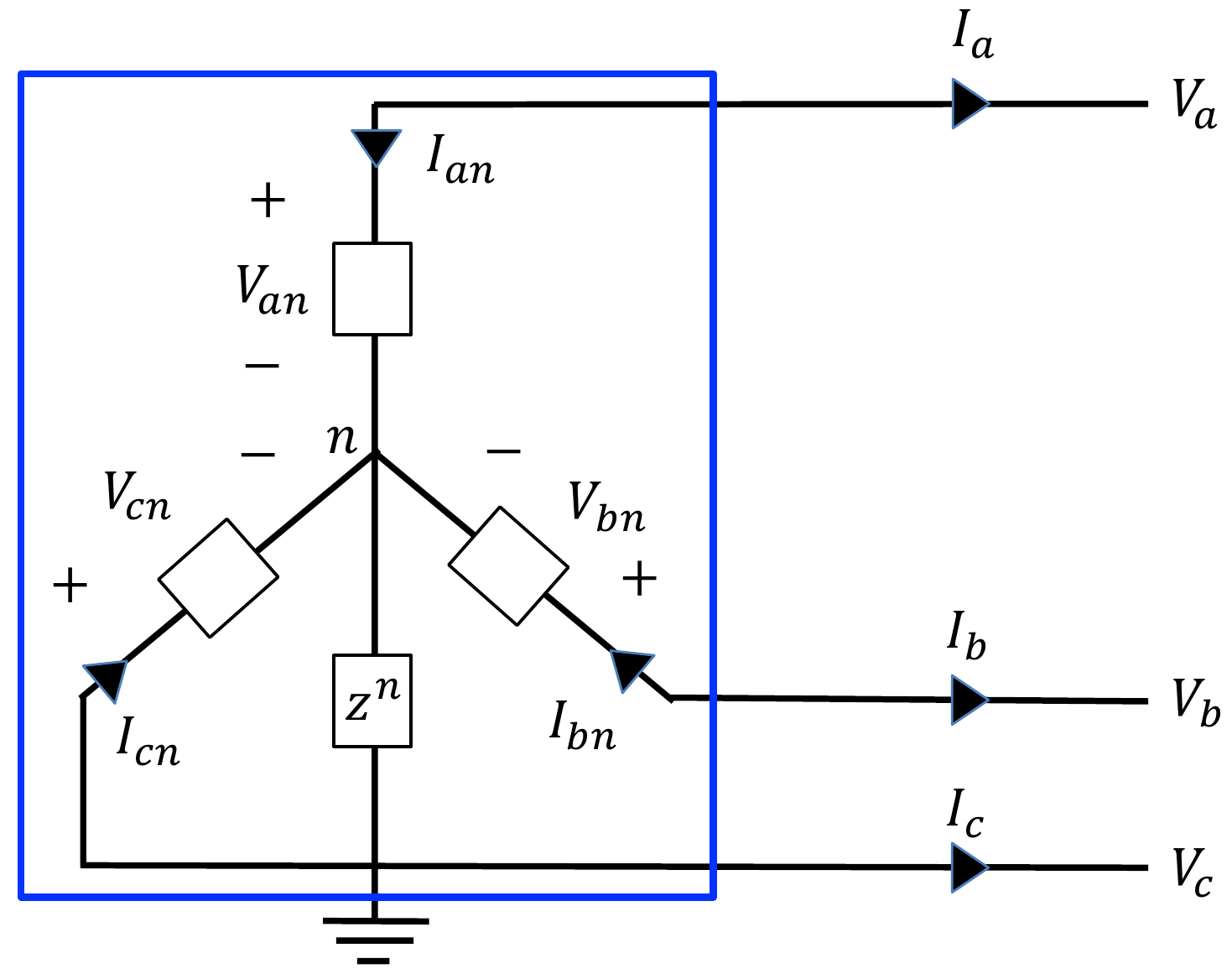}
 }
 \qquad\qquad
\subfigure [$\Delta$ configuration] {
   \includegraphics[width=0.35\textwidth] {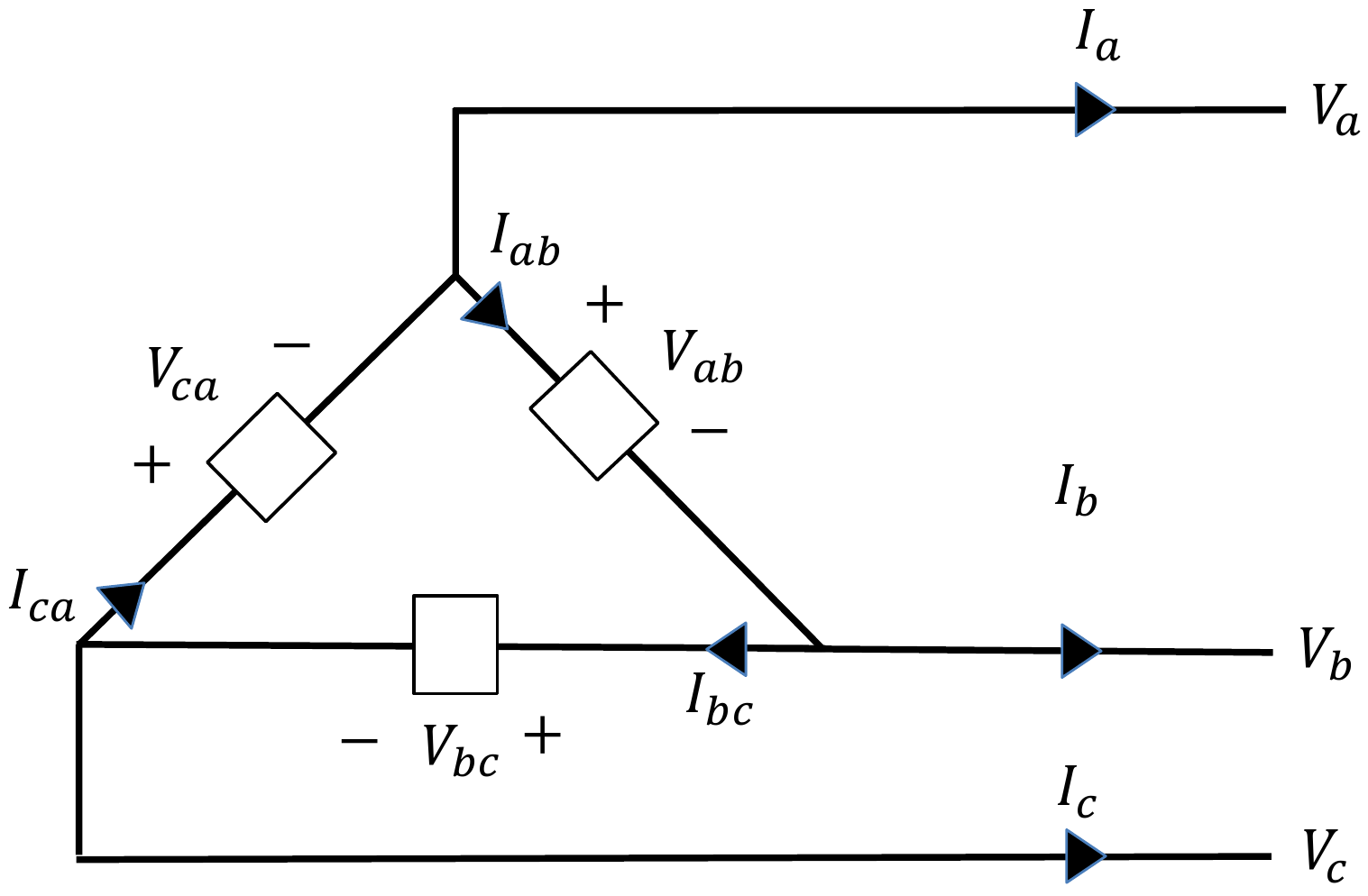}
 }
\caption{Internal and external variables associated with a single-terminal device in $Y$ and $\Delta$ configurations.
	(Even though Figure \ref{ch:mun.1; sec:ComponentModels; fig:YDconfigs.1}(a)
	shows a grounded device, our model allows each $Y$-configured device to be grounded or ungrounded
	and, if grounded, $z^n$ to be zero or nonzero.)}
\label{ch:mun.1; sec:ComponentModels; fig:YDconfigs.1}
\end{figure}
The internal behavior of a single-terminal device shown in 
Figure \ref{ch:mun.1; sec:ComponentModels; fig:YDconfigs.1}
is described in terms of its internal variables:
\bi
\item $V^Y := \left(V^{an}, V^{bn}, V^{cn}\right) \in\mathbb C^3$, $I^Y := \left(I^{an}, I^{bn}, I^{cn}\right) \in\mathbb C^3$, 
	$s^Y := \left(s^{an}, s^{bn}, s^{cn}\right) \in\mathbb C^3$ : line-to-neutral voltages, currents, and power across the single-phase devices 
	in $Y$ configuration.  By definition $s^{a n} := V^{a n} \left( I^{a n}\right)^{\sf H}$ is the power across the phase-$a$
	device, etc. 
	The neutral voltage (with respect to a common reference point) is denoted by $V^n$ and is generally nonzero. 
	A $Y$-configured device may or may not have a neutral line which may or may not be grounded
	(Figure \ref{ch:mun.1; sec:ComponentModels; fig:YDconfigs.1} shows the case where the device is grounded
	through an impedance $z^n$).  When present, the current on the neutral line is denoted by $I^n$ in the direction 
	away from the neutral.

\item $V^\Delta := \left(V^{ab}, V^{bc}, V^{ca}\right) \in\mathbb C^3$, $I^\Delta := \left(I^{ab}, I^{bc}, I^{ca}\right) \in\mathbb C^3$, 
	$s^\Delta := \left(s^{ab}, s^{bc}, s^{ca}\right) \in\mathbb C^3$ : line-to-line voltages, currents, and power across the single-phase devices
	in $\Delta$ configuration.   By definition $s^{ab} := V^{ab} \left( I^{ab}\right)^{\sf H}$ is the power across the phase-$a$
	device, etc.
\ei
Note that the direction of power $s^{an}$ or $s^{ab}$ across a single-phase device is defined in
the \emph{direction of the current through} the device.

\paragraph{Internal models.}

The internal behavior of a three-phase device is described by the relation between the internal variables 
$\left( V^Y, I^Y, s^Y \right)$ or between $\left( V^\Delta, I^\Delta, s^\Delta \right)$.  It depends on the property of the single-phase 
device, but not on their configuration nor the presence of a neutral line:
\bee
\item \emph{Voltage source:} 
	An ideal voltage source fixes the internal voltage $V^{Y/\Delta}$ to be a given constant $E^{Y/\Delta}$.
	
\item \emph{Current source:} 
	An ideal current source fixes the internal current $I^{Y/\Delta}$ to be a given constant $J^{Y/\Delta}$.

\item \emph{Power Source:} An ideal power source fixes the internal power 
	$s^{Y/\Delta}$ to be a given constant $\sigma^{Y/\Delta}$.  

\item \emph{Impedance:} 
	An impedance $z^{Y/\Delta}$ is a $3\times 3$ complex matrix. 
	It fixes the relationship between its internal variables to be:
	\begin{align*}
	V^{Y/\Delta} & \ \ = \ \ z^{Y/\Delta} \, I^{Y/\Delta},  &
	s^{Y/\Delta} & \ \ := \ \ \diag \left( V^{Y/\Delta} \left( I^{Y/\Delta} \right)^{\sf H} \right)
	\end{align*}
\eee
The voltage and current $\left( V^{Y/\Delta}, I^{Y/\Delta} \right)$ of a power 
source $\sigma^{Y/\Delta}$ are related quadratically by
\begin{align*}
\sigma^{Y/\Delta} & \ \ = \ \ \diag \left( V^{Y/\Delta} \left( I^{Y/\Delta} \right)^{\sf H} \right)
\end{align*}
The specification as well as internal and external variables of
a three-phase device are summarized in Table \ref{ch:mun.1; sec:ComponentModels; table:notation}.
\begin{table}[htbp]
	\centering
		\begin{tabular} { l | c | c | c | c }
		\hline\hline
			  &  Voltage & Current & Power & Impedances
		\\
		\hline
		Device specification   & $E^{Y/\Delta}$ & $J^{Y/\Delta}$ & $\sigma^{Y/\Delta}$ & 
				$z^{Y/\Delta}$, $y^{Y/\Delta}$, $z^n$
		\\
		Internal variables  & $V^{Y/\Delta}$ & $I^{Y/\Delta}$ & $s^{Y/\Delta}$ & 
		\\
		External variables   & $V$  & $I$ & $s$ & 
		\\
		\hline\hline
	\end{tabular}	
	\caption{Notation for three-phase single-terminal devices
	$(y^{Y/\Delta} := \left( z^{Y/\Delta} \right)^{-1})$.}
\label{ch:mun.1; sec:ComponentModels; table:notation}
\end{table}
As noted above, the internal model does not depend on the configuration nor the presence of
a neutral line.

\paragraph{Conversion rules.}

The external behavior of the single-terminal device shown in Figure \ref{ch:mun.1; sec:ComponentModels; fig:YDconfigs.1}
is described in terms of its terminal variables:
\bi
\item $V := \left(V^{a}, V^{b}, V^{c}\right) \in\mathbb C^3$, $I := \left(I^{a}, I^{b}, I^{c}\right) \in\mathbb C^3$, 
	 $s := \left(s^{a}, s^{b}, s^{c}\right) \in\mathbb C^3$ : terminal voltages, currents, and power.  
	The terminal voltage $V$ is defined with respect to an arbitrary but common reference point, e.g., the ground.  
	The terminal current $I$ is defined in the direction coming out of the 
	device, i.e., $I$ is defined to be the current injection 
	from the device to the rest of the network when it is connected to a bus bar, regardless of whether it generates or
	consumes power.   	
	By definition $s^{a} := V^{a}\, I^{a \sf H}$ is the power across terminal $a$ and the common
	reference point, etc.
\ei
The external model of a device is the relationship between its terminal variables $(V, I, s)$.  
We now derive the conversion rules 
\eqref{part: NetworkModels; ch:UnbalancedNk; sec:overview; eq:ExtBahevior.1} and
\eqref{ch:mun.1; sec:ComponentModels; subsec:GenLoadD; eq:l2landterminal.s1} 
that maps internal variables $\left( V^{Y/\Delta}, I^{Y/\Delta}, s^{Y/\Delta} \right)$ to external 
variables $(V, I, s)$ for devices in $Y$ and $\Delta$ configurations respectively.   These conversion
rules depend only on the configuration and not on the type of devices.
In Section \ref{sec:3pBFM; subsec:ExternalModels}, we apply these conversion rules to the internal 
model of each device to derive its external model.

\noindent\emph{Conversion in $Y$ configuration.} The terminal voltage, 
current, and power $(V, I, s)$ of a $Y$-configured device are 
related to the internal variables $(V^Y, I^Y, s^Y)$ by:
	\begin{align}
	V & \ = \ V^Y \ + \ V^n \textbf 1,	& 	I & \ = \ - I^Y, &
	s & \ = \ - \left( s^Y \ + \ V^n \overline I^Y \right)
	\label{part: NetworkModels; ch:UnbalancedNk; sec:overview; eq:ExtBahevior.1}
	\end{align}
	where $\overline I$ denotes the componentwise complex conjugate of
	the vector $I^Y\in\mathbb C^3$.
	The negative sign on the current and power conversions
	is due to the definition of $\left( I^Y, s^Y \right)$ as internal current and power delivered to the
	single-phase devices whereas $(I, s)$ is defined as the terminal current 
	and power injections out of the three-phase device.
	
	In general the neutral voltage $V^n$ with respect to a common reference
	point is nonzero whether or
	not there is a neutral line and whether or not the neutral is grounded.
	If the neutral is grounded with zero neutral impedance and voltages are defined 
	with respect to the ground, then $V^n = 0$ and $V = V^Y$ and $s = - s^Y$.
	It is important to explicitly include $V^n$ in a network model because
not every device in a network may be grounded or grounded with zero neutral impedance.

\begin{remark}[Total power]
The total terminal power is
\begin{align*}
\textbf 1^{\sf T} s & \ \ = \ \ - \textbf 1^{\sf T} s^Y \ - \ V^n \left(\textbf 1^T \overline I^Y \right)
\end{align*}
The first term $\textbf 1^{\sf T} s^Y$ is the total power delivered across the single-phase devices. 
The second term $\textbf 1^T I^Y$ is the sum of internal line-to-neutral current.  If 
the neutral is grounded through an impedance then $V^n \left( \textbf 1^T \overline I^Y \right)$ is the power delivered
to the neutral impedance.  If the neutral is ungrounded then $\textbf 1^T I^Y = 0$ by KCL and
the second term $V^n \left( \textbf 1^T \overline I^Y \right) = 0$.
\qed
\end{remark}

\noindent\emph{Conversion in $\Delta$ configuration.} 
The relationship between terminal voltage and current $(V, I)$ and internal voltage and current $\left( V^\Delta, I^\Delta\right)$ is:
\begin{subequations}
\begin{align}
V^\Delta & \ \ = \ \ \Gamma\, V, 		\qquad\quad 		I \ \ = \ \ - \Gamma^{\sf T}\, I^\Delta
\label{ch:mun.1; sec:ComponentModels; subsec:GenLoadD; eq:l2landterminal.VI}
\end{align}
Given appropriate vectors $V^\Delta$ and $I$, solutions $V$ and $I^\Delta$
to \eqref{ch:mun.1; sec:ComponentModels; subsec:GenLoadD; eq:l2landterminal.VI}
is provided by  Lemma 
\ref{part:NetworkModels; ch:mun.1; sec:ComponentModels; subsec:GenLoadD; thm:InverseG}.
\bee
\item Given $V^\Delta$, there is a solution $V$ to 
 \eqref{ch:mun.1; sec:ComponentModels; subsec:GenLoadD; eq:l2landterminal.VI}
 if and only if $V^\Delta$ is orthogonal to $\textbf 1$, i.e., 
\bqn
V^{ab} + V^{bc} + V^{ca} & = & 0
\eqn
which expresses Kirchhoff's voltage law.  In that case, there is a subspace 
of solutions $V$ given by
\bq
V & = & \frac{1}{3}\, \Gamma^{\sf T} V^\Delta \ + \ \gamma\, \textbf 1,
	\quad \gamma\in\mathbb \mathbb C
\label{ch:mun.1; sec:ComponentModels; subsec:GenLoadD; eq:l2landterminal.VI.v}
\eq
This amounts to an arbitrary reference voltage for $V$.  The
quantity $\gamma := \frac{1}{3}\textbf 1^{\sf T} V$ is the (scaled) {zero-sequence}
voltage of $V$.  
In most applications we are given a reference voltage (e.g., $V_0 := \alpha_+$ at the reference
bus 0) which will fix the constant $\gamma$. 

\item Given $I$, there is a solution $I^\Delta$ to 
 \eqref{ch:mun.1; sec:ComponentModels; subsec:GenLoadD; eq:l2landterminal.VI}
 if and only if $I$ is orthogonal to $\textbf 1$, i.e., 
\bqn
I^a + I^b + I^c & = & 0
\eqn
which expresses Kirchhoff's current law.  In that case, there is a subspace of $I^\Delta$ that
satisfy \eqref{ch:mun.1; sec:ComponentModels; subsec:GenLoadD; eq:l2landterminal.VI}, given by
\bq
I^\Delta & = & - \frac{1}{3}\, \Gamma I \ + \ \beta\, \textbf 1,
	\quad \beta\in\mathbb C
\label{ch:mun.1; sec:ComponentModels; subsec:GenLoadD; eq:l2landterminal.VI.i}
\eq
where $\beta$ specifies the amount of loop flow in $I^\Delta$ and does not affect the terminal current $I$
since $\Gamma^{\sf T} \beta\textbf 1 = 0$.  
The quantity $\beta := \frac{1}{3}\textbf 1^{\sf T} I^\Delta$ is the (scaled) {zero-sequence}
current of $I^\Delta$.  
\eee

The terminal power injection from the device is $s := \diag\left( VI^{\sf H} \right)$ and the 
internal power delivered across the single-phase devices in the direction $ab$, $bc$, $ca$ 
is $s^\Delta := \diag\left( V^\Delta I^{\Delta \sf H} \right)$.  
Given internal voltage and current $\left( V^\Delta, I^\Delta \right)$ with $\textbf 1^{\sf T} V^\Delta = 0$,
the terminal power $s$ is
(from \eqref{ch:mun.1; sec:ComponentModels; subsec:GenLoadD; eq:l2landterminal.VI}
\eqref{ch:mun.1; sec:ComponentModels; subsec:GenLoadD; eq:l2landterminal.VI.v}):
\begin{align}
s & 	
\ = \  -  \diag\left( \Gamma^\dag \left(V^\Delta I^{\Delta \sf H} \right) \Gamma \right) + \, \gamma \overline I,	
& \textbf 1^{\sf T} V^\Delta & \ = \ 0    
\label{ch:mun.1; sec:ComponentModels; subsec:GenLoadD; eq:l2landterminal.s1a}
\end{align}
where $\overline I$ is the complex conjugate of the terminal current $I = -\Gamma^{\sf T}I^\Delta$ 
and $\gamma\in\mathbb C$ is 
determined by a reference voltage.   
Conversely, given terminal voltage and current $(V, I)$ with $\textbf 1^{\sf T} I = 0$,
the internal power $s^\Delta$ is
(from \eqref{ch:mun.1; sec:ComponentModels; subsec:GenLoadD; eq:l2landterminal.VI}
\eqref{ch:mun.1; sec:ComponentModels; subsec:GenLoadD; eq:l2landterminal.VI.i}):
\begin{align}
s^\Delta &		
	\ = \-  \diag\left( \Gamma \left( V I^{\sf H} \right) \Gamma^\dag \right) + \, \beta V^\Delta,	&
\textbf 1^{\sf T} I & \ = \ 0
\label{ch:mun.1; sec:ComponentModels; subsec:GenLoadD; eq:l2landterminal.s1b}
\end{align}
where $V^\Delta = \Gamma V$ and $\beta\in\mathbb C$ is determined by the 
zero-sequence current of $I^\Delta$.
\label{ch:mun.1; sec:ComponentModels; subsec:GenLoadD; eq:l2landterminal.s1}
\end{subequations}
\begin{remark}[Total power]
\label{ch:mun.1; sec:ComponentModels; subsec:GenLoadD; remark:powers.1}
\bee
\item Given an internal voltage and current $\left( V^\Delta, I^\Delta \right)$, the terminal power vector $s$
does not depend on the zero-sequence current $\beta := \frac{1}{3}\textbf 1^{\sf T} I^\Delta$ but
does depend on the zero-sequence voltage $\gamma := \frac{1}{3}\textbf 1^{\sf T} V$.
Since $I = -\Gamma^{\sf T} I^\Delta$ and hence $\textbf 1^{\sf T} I = 0$,
the total terminal power however is independent of $\gamma$:
\begin{align*}
\textbf 1^{\sf T} \, s & \ \ = \ \  - \textbf 1^{\sf T} \diag\left( \Gamma^\dag \left(V^\Delta I^{\Delta \sf H} \right) \Gamma \right) 
\end{align*}

\item 
Given a terminal voltage and current $(V, I)$, from 
\eqref{ch:mun.1; sec:ComponentModels; subsec:GenLoadD; eq:l2landterminal.s1b},
the internal power vector $s^\Delta$ depends on zero-sequence current $\beta$.  
Since $V^\Delta  = \Gamma\, V$ and hence $\textbf 1^{\sf T} V^\Delta = 0$, the total internal
power however is independent of the loop flow:
\bqn
\textbf 1^{\sf T} s^\Delta & = & - \textbf 1^{\sf T} \diag\left( \Gamma \left( V I^{\sf H} \right) \Gamma^\dag \right)
\eqn
\eee
\qed
\end{remark}

In summary a complete model of a three-phase device is given by
its internal model specifying the relationship among its
internal variables $\left(V^{Y/\Delta}, I^{Y/\Delta}, s^{Y/\Delta} \right)$
and the conversion rules \eqref{part: NetworkModels; ch:UnbalancedNk; sec:overview; eq:ExtBahevior.1}
and \eqref{ch:mun.1; sec:ComponentModels; subsec:GenLoadD; eq:l2landterminal.s1}
between its internal variables and external variables $(V, I, s)$, 
for $Y$ and $\Delta$ configuration respectively.
This model is required to fully specify a network model
(see below) when the application under study needs to determine or optimize 
some of the internal variables such as the current $I_j^{Y/\Delta}$ or 
power $s_j^{Y/\Delta}$ of each of the single-phase devices connected at a 
bus $j$.

When the application does not require internal variables, we can apply the conversion rules
	\eqref{part: NetworkModels; ch:UnbalancedNk; sec:overview; eq:ExtBahevior.1}
	\eqref{ch:mun.1; sec:ComponentModels; subsec:GenLoadD; eq:l2landterminal.s1}	
to the internal models 
to eliminate the internal variables and obtain a relationship between the external variables $(V, I, s)$
in terms of device parameters, such as $E^{Y/\Delta}$
for an ideal voltage source or $z^{Y/\Delta}$ of an impedance, as we
explain next.

\subsection{Devices: external models}
\label{sec:3pBFM; subsec:ExternalModels}

Since we do not need power sources in this paper, to save space, we omit the
derivation of their external model.

\noindent\emph{$Y$ configuration.}
Application of the conversion rule 
\eqref{part: NetworkModels; ch:UnbalancedNk; sec:overview; eq:ExtBahevior.1}
to the internal models of a voltage source, a current source, and an impedance
yields the external models that relate their terminal variables.
The result is summarized in Table 
\ref{part: NetworkModels; ch:UnbalancedNk; sec:overview; table:IdealDevicesY}.
\begin{table}[htbp]
	\centering
		\begin{tabular} {| l | l | l | }
		\hline\hline
		Device  & \multicolumn{2}{c |}{$Y$ configuration} 
		\\
		\hline
		Voltage source &  $V = E^Y + \gamma\textbf 1$ & 
			$s = \diag\left( E^Y  I^{\sf H}\right) + \gamma\overline I$ 
		\\
		Current source & $\ I = -J^Y$ & $s = -\diag\left( V J^{Y \sf H} \right)$ 
		\\
		Power source & $\diag\left(I^{\sf H}\right) (V - \gamma\textbf 1) = - \sigma$ 
			& $s = - \sigma^Y + \gamma \overline I$  
		\\
		Impedance  & $ V = - z^Y I + \gamma\textbf 1$ & 
			$s = - \diag\left( V \left( V - \gamma \textbf 1 \right)^{\sf H} y^{Y\sf H} \right)$ 
		\\
		\hline\hline
	\end{tabular}	
	\caption{External models of ideal single-terminal devices in $Y$ configuration
	$(\gamma = V^n)$.}
\label{part: NetworkModels; ch:UnbalancedNk; sec:overview; table:IdealDevicesY}
\end{table}
In the table 	
$\gamma = V^n$ is the neutral voltage. These models for
ideal sources do not rely on a common and often implicit assumption that all neutrals are grounded through an impedance $z^n\in\mathbb C$, which may or 
may not be zero, and voltages are defined with respect to the ground.
If this assumption holds then $\gamma = V^n = - z^n \left( \textbf 1^{\sf T}I \right)$ by KCL.  
If, in addition, all neutrals are directly grounded, i.e., $z^n = 0$,
then $\gamma = V^n = 0$ for all $Y$-configured devices.
For a typical three-phase analysis problem, $\gamma$ for all $Y$-configured device 
needs to be specified (see Remark \ref{ch:mun.1; sec:bim; subsec:OverallNk; remark:DeviceSpec}).

\noindent\emph{$\Delta$ configuration.}
The external models of $\Delta$-configured devices can be derived by
applying the conversion rule 
\eqref{ch:mun.1; sec:ComponentModels; subsec:GenLoadD; eq:l2landterminal.s1}	
to their internal models.
\bee
\item \emph{Voltage source $\left(E^\Delta, \gamma \right)$:}
Applying the conversion rule $V^\Delta = \Gamma\, V$ in
\eqref{ch:mun.1; sec:ComponentModels; subsec:GenLoadD; eq:l2landterminal.VI}
to the internal model $V^\Delta = E^\Delta$ of an ideal voltage source, 
we obtain the following external model that relate the terminal voltage, current 
and power $(V, I, s)$:
\begin{subequations}
\begin{align}
V & \ = \ \frac{1}{3} \Gamma^{\sf T} E^\Delta \ + \ \gamma \textbf 1, 	\quad
\textbf 1^{\sf T} I  \ = \ 0
\\
s & \ = \ \frac{1}{3} \diag\left( \Gamma^{\sf T} E^\Delta I^{\sf H} \right) \ + \ \gamma \overline I
\end{align}
provided $\textbf 1^{\sf T} E^\Delta = 0$, 
\label{ch:mun.1; sec:ComponentModels; subsec:GenLoadD; eq:VoltageScr.1a}
\end{subequations}
where $\gamma\in\mathbb C$ is fixed by a given reference voltage.
To specify the external model of an ideal voltage source is to fix the two parameters
$\left( E^\Delta, \gamma\right)$.  Its terminal current and power $(I, s)$ will be 
determined by the interaction of its external model 
\eqref{ch:mun.1; sec:ComponentModels; subsec:GenLoadD; eq:VoltageScr.1a}
with those of other devices on the network through current or power balance
equations.

\item \emph{Current source $J^\Delta$:}
Multiplying $-\Gamma^{\sf T}$ to both sides of the internal model $I^\Delta = J^\Delta$ 
of an ideal current source
and applying the conversion rule $I = - \Gamma^{\sf T} I^\Delta$ in
\eqref{ch:mun.1; sec:ComponentModels; subsec:GenLoadD; eq:l2landterminal.VI}, 
we obtain the external model:
\begin{align}
I & \ = \ - \Gamma^{\sf T}\! J^\Delta,  &
s & \ = \ - \diag \left( VJ^{\Delta \sf H}\Gamma \right)      & & 
\label{ch:mun.1; sec:ComponentModels; subsec:GenLoadD; eq:CurrentScr.1a}
\end{align}
To specify the external model of an ideal current source is to fix the internal
current $J^\Delta$ (which also fixes its zero-sequence current
$\beta := \frac{1}{3}\textbf 1^{\sf T} J^\Delta$).
Its terminal voltage and power $(V, s)$ will be determined
by the interaction of its external model 
\eqref{ch:mun.1; sec:ComponentModels; subsec:GenLoadD; eq:CurrentScr.1a}
with those of other devices on the network through current or power balance equations.

\item \emph{Impedance $z^\Delta$:}
Define the admittance matrix $y^\Delta := \left( z^\Delta \right)^{-1}$.  
Substituting into the internal model
$y^\Delta V^{\Delta} = I^{\Delta}$ of an impedance, multiplying
both sides by $-\Gamma^{\sf T}$ and applying the conversion rule $I = -\Gamma^{\sf T} I^\Delta$, we get
\begin{subequations}
\begin{align}
I & \ \ = \ \ - {Y^\Delta} V		
\label{ch:mun.1; sec:ComponentModels; subsec:GenLoadD; eq:AdmittanceD.1}
\end{align}
where $Y^\Delta$ is a complex symmetric Laplacian matrix given by
\bqn
Y^\Delta \ := \ \Gamma^{\sf T} \, y^\Delta \, \Gamma \ = \
\begin{bmatrix} y^{ca} + y^{ab} & -y^{ab} & -y^{ca}  \\  -y^{ab} & y^{ab}+y^{bc} & -y^{bc}  \\  
	-y^{ca} & -y^{bc} & y^{bc} + y^{ca} \end{bmatrix}
 \label{ch:mun.1; sec:ComponentModels; subsec:GenLoadD; eq:VoltageScr.1b}
\eqn
Note that the terminal current $I$ given by 
\eqref{ch:mun.1; sec:ComponentModels; subsec:GenLoadD; eq:AdmittanceD.1}
satisfies $\textbf 1^{\sf T} I = 0$.
The terminal power injection $s$ can be expressed in terms of $V$:
\begin{align}
s & \ \ = \ \ \diag \left(VI^{\sf H} \right) \ \ = \ \ - \diag\left( V V^{\sf H} Y^{\Delta \sf H} \right)
\end{align}
\label{ch:mun.1; sec:ComponentModels; subsec:GenLoadD; eq:Impedance.1}
\end{subequations}
\eee
The external models 
\eqref{ch:mun.1; sec:ComponentModels; subsec:GenLoadD; eq:VoltageScr.1a}
\eqref{ch:mun.1; sec:ComponentModels; subsec:GenLoadD; eq:CurrentScr.1a}
\eqref{ch:mun.1; sec:ComponentModels; subsec:GenLoadD; eq:Impedance.1}
of ideal $\Delta$-configured devices are
summarized in Table \ref{part: NetworkModels; ch:UnbalancedNk; sec:overview; table:IdealDevicesD}.
\begin{table}[htbp]
	\centering
		\begin{tabular} {| l | l | l | }
		\hline\hline
		Device & \multicolumn{2}{c |} {$\Delta$ configuration}
		\\
		\hline
		Voltage source & 
			$V = \frac{1}{3}\Gamma^{\sf T} E^\Delta + \gamma\textbf 1$, $\textbf 1^{\sf T} I = 0$ & 
			$s = \frac{1}{3}\diag\left( \Gamma^{\sf T} E^\Delta I^{\sf H} \right) + \gamma \overline I$
		\\
		Current source   & 
			$I = -\Gamma^{\sf T} J^\Delta$ & $s = -\diag\left( V J^{\Delta \sf H} \Gamma \right)$
		\\
		Power source &
			$\sigma^\Delta = \diag\left( \Gamma V I^{\Delta \sf H} \right)$ & 
		\\
		Impedance    & $I = - Y^\Delta V$ & 	
		    	$s = - \diag\left( VV^{\sf H} Y^{\Delta \sf H} \right)$
		\\
		\hline\hline
	\end{tabular}	
	\caption{External models of ideal  single-terminal devices in $\Delta$ configuration
	$\left(\gamma := \frac{1}{3} \textbf 1^{\sf T} V, \beta: = \frac{1}{3}\textbf 1^{\sf T} I^\Delta\right)$.}
	\label{part: NetworkModels; ch:UnbalancedNk; sec:overview; table:IdealDevicesD}
\end{table}

\begin{remark}[Non-ideal devices]
For simplicity of exposition, we have presented in this paper the
external models of only ideal devices where the internal series impedances 
of voltage sources and shunt admittances of current sources are assumed zero.  
These models can be extended to non-ideal devices (see \cite{Low2022}).
\qed
\end{remark}

\begin{remark}[$\Delta$-$Y$ transformation]
\label{sec:3pBFM; subsec:ExternalModels; remark:DYtransformation.1}
From the external model  
\eqref{ch:mun.1; sec:ComponentModels; subsec:GenLoadD; eq:VoltageScr.1a} of an ideal
$\Delta$-configured voltage source and that of an $Y$-configured voltage source in 
Table \ref{part: NetworkModels; ch:UnbalancedNk; sec:overview; table:IdealDevicesY}, 
the $Y$ equivalent of $\left( E^\Delta, \gamma \right)$, \emph{not} necessarily balanced, is 
given by
\begin{align*}
E^Y & \ \ := \ \ \frac{1}{3}\Gamma^{\sf T} E^\Delta, 	& 	V^n & \ \ := \ \ \gamma	& & 
\end{align*}
If $E^\Delta$ is balanced then 
$\Gamma^{\sf T}\, E^\Delta = (1-\alpha^2) E^\Delta = {\sqrt 3}\, e^{-\ii \pi/6}\, E^\Delta$ 
and the $Y$ equivalent $E^Y$ reduces to the familiar expression:
\begin{align*}
E^Y & \ \ = \ \  \frac{1}{\sqrt{3}\, e^{\ii{\pi}/6}}\, E^\Delta
\end{align*}

Similarly an {ideal} $\Delta$-configured current source $J^\Delta$ has an $Y$ equivalent $J^Y$ 
given by
\begin{align*}
J^Y & \ \ = \ \ - \Gamma^{\sf T} J^\Delta
\end{align*}
If $J^\Delta$ is balanced then 
\begin{align*}
J^Y & \ \ = \ \  - (1-\alpha^2) J^\Delta \ \ = \ \ - \frac{\sqrt{3}}{e^{\ii{\pi}/6}}\, J^\Delta
\end{align*}
\qed
\end{remark}

\subsection{Three-phase line model}
\label{part: NetworkModels; ch:UnbalancedNk; sec:ComponentModels; subsec:3pLine}

A three-phase line has three wires one for each phase $a, b, c$.  It may also
have a neutral wire which may be grounded at one or both ends if the 
device connected to that end of the line is in $Y$ configuration.
The electromagnetic interactions among the electric charges 
in wires of different phases couple the voltages on and currents in these wires.
The relation between the voltages and currents in these phases can be modeled 
by a linear mapping that depends on the line characteristics.
For simplicity we will restrict ourselves to a three-wire line model that 
takes into account the effect of neutral or earth return on the impedance of a transmission line.
All analysis extends to four-wire models
(including a neutral line) or five-wire models (including a neutral line and 
the ground return) almost without change with proper definitions that include 
neutral and ground variables.

A three-phase line $(j,k)$ is characterized by three $3\times 3$ matrices
$\left( y_{jk}^s, y_{jk}^m, y_{kj}^m \right)$ where $y_{jk}^s$ is the
\emph{series admittance matrix} and $\left( y^m_{jk}, y^m_{kj} \right)$ are
the \emph{shunt admittance matrices}, not necessarily equal.
The terminal voltages $\left( V_j, V_k \right)$ and the {sending-end} currents $\left( I_{jk}, I_{kj} \right)$ 
respectively are related according to
\begin{subequations}
\bq
I_{jk} & = & y^s_{jk}\left( V_j - V_k \right) \ + \ y^m_{jk}\, V_j
\label{ch:umn; sec:ComponentModels; subsec:3pLines; eq:3pmodel.3}
\\
I_{kj} & = & y^s_{jk}\left( V_k - V_j \right) \ + \ y^m_{kj}\, V_k 
\eq
Note that the voltages $(V_j, V_k)$ and currents $(I_{jk}, I_{kj})$
are terminal voltages and currents regardless of whether the three-phase devices connected to
terminals $j$ and $k$ are in $Y$ or $\Delta$ configuration. 

To describe the relationship between the sending-end line power and the voltages
$\left( V_j, V_k \right)$, define the matrices $S_{jk}, S_{kj} \in\mathbb C^{3\times 3}$ by
\begin{align}
S_{jk} &\ \ := \ \ V_j \left( I_{jk} \right)^{\sf H}, 	&
\label{ch:umn; sec:ComponentModels; subsec:3pLines; eq:3pmodel.1a}
S_{kj} & \ \ := \ \ V_k \left( I_{kj} \right)^{\sf H} 	& &
\end{align}
The
\label{ch:umn; sec:ComponentModels; subsec:3pLines; eq:3pmodel.IS}
\end{subequations}
 three-phase sending-end line power from terminals $j$ to $k$ along the line is the vector 
 $\diag\left( S_{jk} \right)$ of diagonal entries and that in the opposite direction is the vector $\diag\left( S_{kj} \right)$.  
The off-diagonal entries of these matrices represent electromagnetic coupling between phases.

\subsection{Network model}
\label{sec:3pBFM; subsec:NetworkModel}

Let $(V, I, s) := \left( V_j, I_j, s_j,  j\in\overline N \right) \in\mathbb C^{3(N+1)}$ be terminal (nodal) variables
over the entire network.  A network equation is a relationship between the {terminal} voltage and
current $(V, I)$ or a relationship between the {terminal} voltage and power $(V, s)$, independent 
of the internal $Y$ or $\Delta$ configurations of the three-phase devices that are connected by the lines.
In both cases the extension of the line model \eqref{ch:umn; sec:ComponentModels; subsec:3pLines; eq:3pmodel.IS}
to a network is simply the nodal current or power balance equations:
\begin{align*}
I_j  & \ \ = \ \ \sum_{k: j\sim k} I_{jk}, 	& s_j  & \ \ = \ \  \sum_{k: j\sim k} \diag\left( S_{jk} \right), 	&
j & \ \ \in \ \ \overline N
\end{align*}
where $S_{jk}$ are matrices defined in 
\eqref{ch:umn; sec:ComponentModels; subsec:3pLines; eq:3pmodel.1a}.
In this paper we focuses on the current balance equation which, 
using \eqref{ch:umn; sec:ComponentModels; subsec:3pLines; eq:3pmodel.3}, is:
\begin{subequations}
\bq
I_j & = & \sum_{k: j\sim k} \left( y^s_{jk} + {y_{jk}^m} \right) V_j - \sum_{k: j\sim k} y^s_{jk} V_k,
	\quad j\in\overline N
\label{ch:mun.1; sec:bim; subsec:I-V relationship; eq:I=YV.1a}
\eq
Note 
that $I_j$ is the net current injection.\footnote{If there is a nodal shunt admittance load $y_j^\text{sh}$, 
e.g., a capacitor bank, in addition to a device whose terminal injection is $\tilde I_j$, then the net injection 
from bus $j$ to the rest of the network is $I_j = \tilde I_j - y_j^\text{sh} V_j$.  This assumes that $y_j^\text{sh}$ 
connects bus $j$ to the ground and the terminal voltage $V_j$ is defined with respect to the ground.}
In vector form, this relates the bus current vector $I := (I_0, \dots, I_N)$ to the bus voltage vector $V := (V_0, \dots, V_N)$:
\bq
I & = & YV
\label{ch:mun.1; sec:bim; subsec:I-V relationship; eq:I=YV.2}
\eq
in terms of a $3(N+1) \times 3(N+1)$ admittance matrix $Y$ where its $3 \times 3$ submatrices are given by
\bq
Y_{jk} & = & \left\{ \begin{array}{lcl}
			- y^s_{jk},  & & j\sim k \ \ (j\neq k)  \\
			\sum_{l: j\sim l}\, y^s_{jl} \ + \  y_{jj}^m,   & & j=k   \\
			0   & & \text{otherwise}
			\end{array}  \right.
\label{eq:defYjk}
\eq
\label{eq:defY}
\end{subequations}

An overall network model consists of (see Figure \ref{fig-IREP-Overall3pModelingSteps}):
\bee
\item A network model \eqref{eq:defY} that relates terminal voltage and
    current $(V, I)$.  

\item A device model for each three-phase device $j$.  
    This can either be:
    \bi
    \item An internal model together with the conversion rules 
    \eqref{part: NetworkModels; ch:UnbalancedNk; sec:overview; eq:ExtBahevior.1}\eqref{ch:mun.1; sec:ComponentModels; subsec:GenLoadD; eq:l2landterminal.s1}
    in Section \ref{sec:3pBFM; subsec:InternalModels}; or 
    \item An external model summarized in Tables 
    \ref{part: NetworkModels; ch:UnbalancedNk; sec:overview; table:IdealDevicesY} and 
	\ref{part: NetworkModels; ch:UnbalancedNk; sec:overview; table:IdealDevicesD} 
    in Section \ref{sec:3pBFM; subsec:ExternalModels} 	
    when only terminal quantities are needed.
    \ei
 \eee

\section{Three-phase analysis}
\label{sec:3pAnalysis}

We now formulate a general three-phase analysis problem using the
overall model of Section \ref{sec:3pModel}.
Consider a three-phase network $G := (\overline N, E)$ where each line $(j,k)\in E$ is characterized by 
$3\times 3$ series and shunt admittance matrices $\left( y_{jk}^s, y_{jk}^m, y_{kj}^m\right)$.   
At each bus $j\in\overline N$ we assume, without loss of generality, there is a single
 three-wire device in either $Y$ or $\Delta$ configuration.

\noindent\emph{Three-phase devices.}
Partition $\overline N$ into 6 disjoint subsets:
\bi
\item $N_v^{Y/\Delta}$: buses with ideal voltage sources in $Y$ or $\Delta$ configurations.
\item $N_c^{Y/\Delta}$: buses with ideal current sources in $Y$ or $\Delta$ configurations.
\item $N_i^{Y/\Delta}$: buses with impedances in $Y$ or $\Delta$ configurations.
\ei
\begin{table*}[htbp]
	\centering
		\begin{tabular} {| c | c | l | c | }
		\hline\hline
		Buses $j$  & Specification & {External model} & Unknowns
		\\
		\hline
		$N_v^Y$  & ${\color{blue} {V_j^Y := E_j^Y}}$, ${\color{blue}{\gamma_j}}$ &  
			 {\color{blue}{$V_j \ = \ E_j^Y + \gamma_j\textbf 1$}} 
			 & $\left( I_j, I_j^Y \right)$
		\\
		$N_v^\Delta$  & ${\color{blue} {V_j^\Delta := E_j^\Delta}}$,  {\color{blue}{$\gamma_j, \beta_j$}},   &  
			 {\color{blue}{$V_j \ = \ \Gamma^\dag E_j^\Delta + \gamma_j\textbf 1$}}	
			 & $\left( I_j, I_j^\Delta \right)$
		\\
		$N_c^Y$  & ${\color{blue}{I_j^Y := J_j^Y, \gamma_j}}$ 
			& {\color{blue}{$I_j  \ = \ - J_j^Y$}}
			& $\left( V_j, V_j^Y \right)$
	\\
		$N_c^\Delta$  & ${\color{blue}{I_j^\Delta := J_j^\Delta}}$   
			&  {\color{blue}{$I_j \ = \ - \Gamma^{\sf T} J_j^\Delta $}}
			& $\left( V_j, V_j^\Delta \right), \left( \gamma_j, \beta_j \right)$
	\\
		$N_i^Y$  & ${\color{blue}{z_j^Y, \ \gamma_j}}$   
			& $I_j \ = \ - {\color{blue}{y_j^Y}} \left( V_j - {\color{blue}{\gamma_j}}\textbf 1\right)$
			& $\left( V_j, I_j \right), \left( V_j^Y, I_j^Y \right)$
	\\
		$N_i^\Delta$  & {\color{blue}{$z_j^\Delta, \ \beta_j$}}
			& $I_j  \ = \ -  {\color{blue}{Y_j^\Delta}} V_j$
			& $\left( V_j, I_j \right), \left( V_j^\Delta, I_j^\Delta \right), \gamma_j$
	\\	
	\hline\hline
	\end{tabular}	
	\caption{Internal and external models of three-phase sources and impedances
	from Tables \ref{part: NetworkModels; ch:UnbalancedNk; sec:overview; table:IdealDevicesY} and
	\ref{part: NetworkModels; ch:UnbalancedNk; sec:overview; table:IdealDevicesD}.
	The three-phase analysis problem is: given the specification in {\color{blue}blue}, compute the
	remaining unknowns in black.
	}
\label{ch:mun.1; sec:bim; subsec:OverallNk; table:devices}
\end{table*}	

\subsection{Device specification}

Associated with each device $j$ are the internal variables 
$\left( V_j^{Y/\Delta}, I_j^{Y/\Delta}, s_j^{Y/\Delta} \right)\in\mathbb C^9$, the terminal variables 
$\left( V_j, I_j, s_j\right)\in\mathbb C^9$, and the variables $\left(\gamma_j, \beta_j\right)\in\mathbb C^2$.
Some of these variables are specified in a three-phase analysis problem and the others are computed
from network equations and conversion rules.  We now describe which of these variables are 
specified for each device in a typical three-analysis problem (without constant-power devices).    The result is summarized in 
Table \ref{ch:mun.1; sec:bim; subsec:OverallNk; table:devices}.  These requirements may need 
to be modified depending on the details of a problem. 
\bee
\item \emph{Voltage source $j$}: It is specified by its internal voltage and a 
    parameter $\left(E_j^{Y/\Delta}, \gamma_j \right)$ where 
    $\gamma_j := V_j^n$ is the neutral voltage if $j$ is in $Y$ configuration and
	$\gamma_j := \frac{1}{3}\textbf 1^{\sf T} V_j$
	is the zero-sequence component of the terminal voltage if $j$ is in $\Delta$
	configuration. For a $\Delta$-configured voltage source, the zero-sequence
	current $\beta_j$ 
	also needs to be specified in order to determine the internal current $I_j^\Delta$
	from the terminal current $I_j$.
	
\item \emph{Current source $j$}: It is specified by its internal current 
    $J_j^{Y/\Delta}$.  For a $Y$-configured current source, its neutral
    voltage $\gamma_j$ is also specified.  

\item \emph{Impedance $j$}: A $Y$-configured impedance $j$ is specified by 
    its internal impedance $z_j^{Y}$ and the neutral voltage $\gamma_j := V_j^n$.
    A $\Delta$-configured impedance $j$ is specified by $z_j^\Delta$ and its
    zero-sequence current $\beta_j$.
\eee
The specification of each device also comes with an external model that
relates its terminal variables in terms of the specified parameters, as
shown in Table \ref{ch:mun.1; sec:bim; subsec:OverallNk; table:devices}.

\subsection{Analysis problem}
\label{sec:3pAnalysis; subsec:Problem}

A three-phase analysis problem is: given devices specified as in 
Table \ref{ch:mun.1; sec:bim; subsec:OverallNk; table:devices} connected
by three-phase lines with given admittance matrices 
$\left( y_{jk}^s, y_{jk}^m, y_{jk}^m \right)$, compute the remaining 
unknowns for each bus $j$ listed in the last column of the table.
The general solution strategy is to use the external models in 
Table \ref{ch:mun.1; sec:bim; subsec:OverallNk; table:devices} and the
network equation $I = YV$ in \eqref{eq:defY}
to compute terminal voltages and currents $\left( V_j, I_j \right)$.  
Internal variables $\left( V_j^{Y/\Delta}, I_j^{Y/\Delta} \right)$ as well 
as $\left(\gamma_j, \beta_j\right)$ can then be determined by the conversion 
rules.

Specifically let $N_v := N_v^Y \cup V_v^\Delta$, $N_c := N_c^Y \cup V_c^\Delta$, and 
$N_i := N_i^Y \cup V_i^\Delta$ be the set of buses with, respectively, voltage sources, current
sources, and impedances.  
With a slight abuse of notation define the following (column) vectors of terminal voltages and currents:
\begin{align*}
\left( V_v, I_v \right) & \ \ := \ \ \left( V_j, I_j, \,  j\in N_v \right)	\\
\left( V_c, I_c \right)& \ \ := \ \ \left( V_j, I_j, \, j\in N_c \right)	\\
\left( V_i, I_i \right) & \ \ := \ \ \left( V_j, I_j, \, j\in N_i \right)
\end{align*}
Then $I = YV$ becomes
\begin{align}
\begin{bmatrix} I_v \\ I_c  \\  I_i  \end{bmatrix} & \ \ = \ \ 
\underbrace{ \begin{bmatrix} Y_{vv} & Y_{vc} & Y_{vi}  \\  Y_{cv} & Y_{cc} & Y_{ci}  \\  Y_{iv} & Y_{ic} & Y_{ii}  
	\end{bmatrix} }_{Y}
\begin{bmatrix} V_v \\ V_c  \\  V_i  \end{bmatrix}
\label{ch:mun.1; sec:bim; subsec:OverallNk; eq:I=YV.1}
\end{align}
where the admittance matrix $Y$ is defined in \eqref{eq:defY}.  
The three-phase network analysis problem is then:
\bee
\item Given the specification of voltage sources, current sources and impedances
	in Table \ref{ch:mun.1; sec:bim; subsec:OverallNk; table:devices}, solve
	\eqref{ch:mun.1; sec:bim; subsec:OverallNk; eq:I=YV.1} for the terminal voltage 
	$V_{-v} := \left( V_c, V_i \right)$ and current $I_{-c} := \left(I_v, I_i \right)$. 

\item From the terminal voltage $V := \left( V_v, V_c, V_i \right)$ and current $I := \left( I_v, I_c, I_i \right)$,
	the internal voltages, currents and powers of $Y$-configured devices are obtained from the
	conversion rule \eqref{part: NetworkModels; ch:UnbalancedNk; sec:overview; eq:ExtBahevior.1} 
	$\left(\gamma_j = V_j^n\right)$:
	\begin{align*}
	V_j^Y & \ \ = \ \ V_j \ - \ \gamma_j \textbf 1, 		& 	I_j^Y & \ = \ - I_j	\\
	s_j^Y & \ \ = \ \ \diag\left( V_j I_j^{\sf H} \right), 	&	j & \ \in \  N_v^Y \cup N_c^Y \cup N_i^Y
	\end{align*}

\item Those of $\Delta$-configured devices can be computed by applying the conversion rule 
\eqref{ch:mun.1; sec:ComponentModels; subsec:GenLoadD; eq:l2landterminal.s1} to $(V, I)$:
	\begin{align*}
	V_j^\Delta & \ \ = \ \ \Gamma V_j, 		& 	
	I_j^\Delta & \ = \ - \Gamma^{\sf T\dag} I_j \ + \ \beta_j \textbf 1	\\
	s_j^\Delta & \ \ = \ \ \diag\left( V_j^\Delta I_j^{\Delta \sf H} \right), 	&
		j & \ \in \  N_v^\Delta \cup N_c^\Delta \cup N_i^\Delta
	\end{align*}
\eee
The main task is to solve the network equation 
\eqref{ch:mun.1; sec:bim; subsec:OverallNk; eq:I=YV.1} in Step 1 for terminal 
voltage and current $(V, I)$ (see \cite{Low2022} for more details).
This generally can only be solved numerically.

The result of the analysis determines both internal and terminal variables
$\left( V_j^{Y/\Delta}, I_j^{Y/\Delta}, s_j^{Y/\Delta} \right)$, 
$\left( V_j, I_j, s_j\right)$ and $\left(\gamma_j, \beta_j\right)$ at every
bus $j$.  We make a few remarks.
\begin{remark}[Voltage $\gamma_j$]
\label{ch:mun.1; sec:bim; subsec:OverallNk; remark:DeviceSpec}
\bee
\item \emph{Parameter $\gamma_j$ for $Y$-configured devices.} 
    The voltage parameter $\gamma_j$ needs to be specified for every $Y$-configured 
    device. By that, we mean information additional 
    to the models in Table \ref{ch:mun.1; sec:bim; subsec:OverallNk; table:devices} 
    is available to determine the value of $\gamma_j$ for that device.
    It may be specified directly, or more likely, indirectly.  
    For instance if the neutral of a $Y$-configured device is grounded and all 
    voltages are defined with respect to the ground, then 
    $\gamma_j = V_j^n = -z_j^n\left( \textbf 1^{\sf T} I_j \right)$, which
    allows the elimination of $\gamma_j$ from the model.  If the neutral is 
    grounded directly (i.e., $z_j^n=0$), then $\gamma_j = 0$.  
    If the neutral is not grounded but the internal 
	voltage $V_j^Y$ is known to be balanced, i.e., $\textbf 1^{\sf T} V_j^Y = 0$, 
	then $\gamma_j := \frac{1}{3}\textbf 1^{\sf T} V_j$.  
    For a $Y$-configured current source, $\gamma_j$ is usually not needed
    to determine its terminal voltage $V_j$, but needed to compute its internal
    voltage $V_j^Y = V_j - \gamma_j\textbf 1$ from the terminal voltage $V_j$.

\item \emph{Variable $\gamma_j$ for $\Delta$-configured devices.} 
    For a $\Delta$-configured voltage source, the zero-sequence voltage 
    $\gamma_j := \frac{1}{3}\textbf 1^{\sf T}V_j$ needs to be specified,
    e.g., by specifying one of its terminal voltages, say, $V_j^a$.  
    For a $\Delta$-configured current source or impedance, $\gamma_j$
    can be determined once its terminal voltage $V_j$ is determined from the
    network equation $I = YV$.

\item \emph{Neutral voltage $\gamma_j$ and zero-sequence voltage.}	For any
    $Y$-configured device, we have 
	\begin{align*}
	V_j & \ \ = \ \ V_j^Y \ + \ V_j^n \textbf 1
	\end{align*}
	The parameter $\gamma_j := V_j^n$ may or may not equal the
	zero-sequence voltage $\frac{1}{3}\textbf 1^{\sf T} V_j$.  They are equal if and only
	if the internal voltages have no zero-sequence component since
	$\frac{1}{3} \textbf 1^{\sf T} V_j =  \frac{1}{3} \textbf 1^{\sf T} V_j^Y + V_j^n$.
\eee
\qed
\end{remark}

\section{Balanced network}
\label{part:networks; ch:mun.bim; sec:BalancedNk}

In this section we show that, if the voltage sources, current sources, and impedances are 
balanced 
and the lines are decoupled, then the three-phase network is equivalent to a per-phase network and
the analysis problem in Section \ref{sec:3pAnalysis} can be solved by analyzing the
simpler per-phase network.  

The intuition is as follows.
In a balanced three-phase network, positive-sequence voltages and 
currents are in span$(\alpha_+)$ and $\alpha_+$ 
is an eigenvector of $\Gamma$ and $\Gamma^{\sf T}$.  This means that the 
transformation of balanced voltages and currents under $\Gamma, \Gamma^{\sf T}$ reduces to a scaling of 
these variables by their eigenvalues $1-\alpha$ and $1-\alpha^2$ respectively.  
The voltage and current at every point in a network can be written
as linear combinations of transformed source voltages and source currents,
transformed by $\left( \Gamma, \Gamma^{\sf T} \right)$ and line admittance matrices.
Therefore if the source voltages and source currents are balanced
positive-sequence sets and lines are identical and phase-decoupled, then
the transformed voltages and currents remain in span$(\alpha_+)$ and hence 
are balanced positive-sequence sets.

In Section \ref{part:networks; ch:mun.bim; sec:BalancedNk; subsec:3pAnalysis}
we describe how the balanced nature of voltage and current sources simplifies
the three-phase analysis problem formulated in Section \ref{sec:3pAnalysis}.
In Section \ref{part:networks; ch:mun.bim; sec:BalancedNk; subsec:1pNetwork}
we describe the positive-sequence per-phase network.  
In Section \ref{part:networks; ch:mun.bim; sec:BalancedNk; subsec:1pAnalysis}
we describe per-phase analysis under the assumption that the neutral voltages
$V_j^n$ of all $Y$-configured devices are zero and the zero-sequence voltages
$\gamma_j$ of all $\Delta$-configured voltage sources are zero, and justify the
procedure in Theorem \ref{ch:mun.1; sec:bim; subsec:BalancedNk; thm:1phiNk.C1C2}.
In Section \ref{part:networks; ch:mun.bim; sec:BalancedNk; subsec:extension}
we extend the per-phase analysis and Theorem \ref{ch:mun.1; sec:bim; subsec:BalancedNk; thm:1phiNk.C1C2}
to the case without this assumption.
In Section \ref{part:networks; ch:mun.bim; sec:BalancedNk; subsec:equivalence}
we prove Theorem \ref{ch:mun.1; sec:bim; subsec:BalancedNk; thm:1phiNk.C1C2}.

\subsection{Problem formulation}
\label{part:networks; ch:mun.bim; sec:BalancedNk; subsec:3pAnalysis}

\noindent\textbf{Balanced devices.}
Three-phase devices are balanced positive-sequence sets if the voltage and current sources
are in span$(\alpha_+)$ and impedances are balanced (identical) across phases.  Then their 
internal models in Table \ref{ch:mun.1; sec:bim; subsec:OverallNk; table:devices} reduce to 
those specified in Table \ref{ch:mun.1; sec:bim; subsec:BalancedNk; table:buses} with parameters
$\lambda_j, \mu_j, \epsilon_j\in\mathbb C$.  
The external models in Table \ref{ch:mun.1; sec:bim; subsec:BalancedNk; table:buses} 
are obtained by substituting these specifications into the external models in
Table \ref{ch:mun.1; sec:bim; subsec:OverallNk; table:devices} and applying 
(Theorem \ref{part:NetworkModels; ch:mun.1; sec:ComponentModels; subsec:GenLoadD; thm:InverseG})
\begin{align*}
\Gamma \alpha_+ & \ \ = \ \ (1-\alpha) \alpha_+, 	& 	\Gamma^{\sf T} \alpha_+ & \ \ = \ \ \left( 1 - \alpha^2 \right) \alpha_+
\\
\Gamma^\dag & \ \ = \ \ \frac{1}{3} \Gamma^{\sf T}, 	& 	\Gamma^{\sf T\dag} & \ \ = \ \ \frac{1}{3} \Gamma
\end{align*}
For example the external model of a $\Delta$-configured impedance is
$I_j = - Y^\Delta V_j$ where the effective matrix 
$Y_j^\Delta := \Gamma^{\sf T} y_j^\Delta\, \Gamma$. 
Since a balanced impedance is $z^\Delta_j = \epsilon_j^{-1}\mathbb I$, we have
\begin{align*}
Y_j^\Delta & \ \ = \ \ \Gamma^{\sf T} y_j^\Delta\, \Gamma \ \ = \ \ 
	\epsilon_j \left(3  \mathbb I - \textbf 1\textbf 1^{\sf T} \right)
\end{align*}
so that, since $\textbf 1^{\sf T} V_j =: 3\gamma_j$, the 
external models of an impedance in $\Delta$ configuration reduces to:
\begin{align*}
I_j & \ \ = \ \ - Y^\Delta V_j  
\ \ = \ \ - 3\epsilon_j\left( V_j \ - \ \gamma_j \textbf 1 \right)
\end{align*}
\begin{table*}[htbp]
	\centering
		\begin{tabular} {| c | l | l | c | l | }
		\hline\hline
		Buses $j$  & \ \ Specification & External model & Vars & Internal vars
		\\
		\hline
		$N_v^Y$  & {\color{blue}{$V_j^Y = \lambda_j\alpha_+$, $\ \gamma_j$}} &  
			${\color{blue}{V_j \ = \ \lambda_j\alpha_+ + \gamma_j\textbf 1}}$ 
			& $I_j$ 	& $I_j^Y = - I_j$
		\\
		$N_v^\Delta$  & {\color{blue} {$V_j^\Delta = \lambda_j\alpha_+$,  $\gamma_j$, $\beta_j$ }} &  
			${\color{blue}{V_j \ = \		
				\frac{1}{3}(1-\alpha^2) \lambda_j\alpha_+ + \gamma_j\textbf 1}}$ 
			& $I_j$	& $I_j^\Delta \ = \ - \Gamma^{\sf T \dag} I_j + \beta_j \textbf 1$
		\\
		$N_c^Y$  & {\color{blue} {$I_j^Y \ = \ \mu_j \alpha_+, \gamma_j$ }}   
			& {\color{blue}{$I_j  \ = \ - \mu_j \alpha_+$}}
			& $V_j$	& $V_j^Y \ = \ V_j - \gamma_j \textbf 1$ 
		\\
		$N_c^\Delta$  & {\color{blue} {$I_j^\Delta \ = \ \mu_j \alpha_+$ }}    
			&  {\color{blue}{$I_j \ = \  - (1-\alpha^2) \mu_j \alpha_+$}}
			& $V_j$	& $V_j^\Delta \ = \ \Gamma V_j, \ \gamma_j := \frac{1}{3} \textbf 1^{\sf T} V_j$ 
		\\	& & & 
			& $\beta_j := \frac{1}{3}\textbf 1^{\sf T} I_j^\Delta$
		\\
		$N_i^Y$  & {\color{blue} {$z_j^Y \ = \ \epsilon_j^{-1} \mathbb I$, $\ \gamma_j$}} &  
			$I_j \ = \ - {\color{blue}{ \epsilon_j}} \left( V_j - {\color{blue}{\gamma_j}} \textbf 1\right)$
			& $\left(V_j, I_j \right)$	& $V_j^Y \ = \ V_j - \gamma_j \textbf 1$, \ $I_j^Y = - I_j$
		\\
		$N_i^\Delta$  & {\color{blue} {$z_j^\Delta \ = \ \epsilon_j^{-1} \mathbb I$, $\ \beta_j$  }}  
			& $I_j = - {\color{blue}{3\epsilon_j}} \left( V_j \ - \ \gamma_j  \textbf 1 \right)$
			& $\left(V_j, I_j \right)$	& $V_j^\Delta \ = \ \Gamma V_j$, \ $\gamma_j := \frac{1}{3} \textbf 1^{\sf T} V_j$
		\\	& & & 
			& $I_j^\Delta \ = \ - \Gamma^{\sf T \dag} I_j + \beta_j \textbf 1$
		\\
		\hline\hline
	\end{tabular}	
	\caption{Internal and external models of balanced positive-sequence sources and impedances.
	}
\label{ch:mun.1; sec:bim; subsec:BalancedNk; table:buses}
\end{table*}	

\noindent\textbf{Balanced admittance matrix $Y$.}
We assume all lines are balanced, i.e., 
\begin{subequations}
\begin{align}
y_{jk}^s & \ \ = \ \ \eta_{jk}^s \mathbb I, 	& 
y_{jk}^m & \ \ = \ \ \eta_{jk}^m \mathbb I, 	& 	y_{kj}^m & \ \ = \ \ \eta_{kj}^m \mathbb I
\label{ch:mun.1; sec:bim; subsec:BalancedNk; eq:1a}
\end{align}
for some constants $\eta_{jk}^s, \eta_{jk}^m, \eta_{kj}^m\in\mathbb C$.
The terminal voltages and currents $V := (V_0, \dots, V_N)$ and $I := (I_0, \dots, I_N)$ are described by \eqref{eq:defY}
which, with balanced lines, reduces to
\bq
I_j & = & 
\sum_{k: j\sim k} \eta_{jk} V_j \ - \ \sum_{k: j\sim k} \eta^s_{jk} V_k,
	\quad j\in\overline N	\quad\quad
\label{ch:mun.1; sec:bim; subsec:BalancedNk; eq:1b}
\eq
where
\label{ch:mun.1; sec:bim; subsec:BalancedNk; eq:1}
\end{subequations}
$\eta_{jk} := \eta_{jk}^s + \eta_{jk}^m$ and $V_j, I_j \in \mathbb C^3$.
This in vector form is $I = YV$.  
%
Define the $(\overline N + 1)\times (\overline N+1)$ \emph{per-phase admittance matrix}
$Y^{1\phi}$ by
\begin{subequations}
\begin{align}
Y^{1\phi}_{jk}  & \ \ := \ \ \left\{ \begin{array}{lcl}
	- \eta^s_{jk}, &  & (j,k) \in E, \ \ (j\neq k)	\\
	\sum_{k: j\sim k} \left( \eta^s_{jk} \ + \ \eta_{jk}^m \right),  & & j = k 	\\
	0  & & \text{otherwise}
	\end{array} \right.
\label{eq:defYjk1phi}
\end{align}
Substituting
\eqref{ch:mun.1; sec:bim; subsec:BalancedNk; eq:1a} into the admittance matrix $Y$ in \eqref{eq:defY} 
for the three-phase network, we can write $Y$ in terms of the per-phase admittance matrix $Y^{1\phi}$ 
using the Kronecker product:
\begin{align}
Y & \ \ = \ \ Y^{1\phi} \otimes \mathbb I
\label{eq:defY1phi}
\end{align}
The relationship $I = YV$ for the three-phase network becomes
\begin{align}
I & \ \ = \ \ \left(Y^{1\phi} \otimes \mathbb I\right) V
\label{eq:I=Y1phiV.1b}
\end{align}
\label{part:networks; ch:mun.bim; sec:BalancedNk; subsec:3pAnalysis; eq:I=Y1phiV.1}
\end{subequations}

\noindent\textbf{Three-phase analysis problem.}
The analysis problem in Section \ref{sec:3pAnalysis} reduces to the
following problem.  
To simplify notation define
\begin{align*}
\hat\alpha_j & \ \ := \ \ \left\{ \begin{array}{lcl}
		1 & \text{ if } & j \in N_v^Y \cup N_c^Y \cup N_i^Y 	\\ 	
		(1-\alpha^2)/3  &  \text{ if } &  j \in N_v^\Delta \qquad (\text{voltage sources})	\\
		(1-\alpha^2)	& \text{ if } & j \in N_c^\Delta	 \qquad (\text{current sources})		\\
		3	& \text{ if } & j \in N_i^\Delta	 \qquad (\text{admittance})		\\
		\end{array} \right.
\end{align*}
Given the following balanced voltage and current sources $\left( V_v, I_c \right)$ and impedance model 
(from Table \ref{ch:mun.1; sec:bim; subsec:BalancedNk; table:buses}):
\begin{subequations}
\begin{align}
V_v & \ \ = \ \ \left( \hat\alpha_j \lambda_j\alpha_+ + \gamma_j\textbf 1, j\in N_v \right)
\label{ch:mun.1; sec:bim; subsec:BalancedNk; eq:3pspec.1a}
\\
I_c & \ \ = \ \ \left( - \hat\alpha_j \mu_j \alpha_+, j\in N_c \right)
\label{ch:mun.1; sec:bim; subsec:BalancedNk; eq:3pspec.1b}
\\
I_j & \ \ = \ \ - \hat\alpha_j \epsilon_j \left( V_j - \gamma_j \textbf 1 \right), \ \ j\in N_i
\label{ch:mun.1; sec:bim; subsec:BalancedNk; eq:3pspec.1c}
\end{align}
our objective
\label{ch:mun.1; sec:bim; subsec:BalancedNk; eq:3pspec.1}
\end{subequations}
is to solve \eqref{part:networks; ch:mun.bim; sec:BalancedNk; subsec:3pAnalysis; eq:I=Y1phiV.1}
for the terminal voltage $V_{-v} := \left( V_c, V_i \right)$ and current $I_{-c} := \left(I_v, I_i \right)$
and then calculate internal voltages and currents as well as $\left( \gamma_j, \beta_j \right)$.

The problem can be solved by substituting 
\eqref{ch:mun.1; sec:bim; subsec:BalancedNk; eq:3pspec.1} into 
\eqref{part:networks; ch:mun.bim; sec:BalancedNk; subsec:3pAnalysis; eq:I=Y1phiV.1}
and computing numerically $\left( V_{-v}, I_{-c} \right)$.
This is Step 1 of the solution procedure in Section 
\ref{sec:3pAnalysis; subsec:Problem}.  Steps 2 and 3 will compute the internal
variables given the terminal variables $(V, I)$.

\begin{remark}[$\Delta$-$Y$ transformation]
The specification \eqref{ch:mun.1; sec:bim; subsec:BalancedNk; eq:3pspec.1} corresponds
to the step of converting all $\Delta$ configured devices to their  $Y$ equivalents.
It generalizes the standard practice of assuming $\gamma_j = 0$ to the case where $\gamma_j$
may be nonzero, because some $Y$-configured devices on the network are not grounded, some 
are grounded through nonzero earthing impedances, and some $\Delta$-configured devices have
nonzero zero-sequence voltages 
(cf. Remark \ref{sec:3pBFM; subsec:ExternalModels; remark:DYtransformation.1}).
\qed
\end{remark}

\subsection{Per-phase network}
\label{part:networks; ch:mun.bim; sec:BalancedNk; subsec:1pNetwork}

We now formalize the alternative solution that solves 
\eqref{part:networks; ch:mun.bim; sec:BalancedNk; subsec:3pAnalysis; eq:I=Y1phiV.1}\eqref{ch:mun.1; sec:bim; subsec:BalancedNk; eq:3pspec.1} 
using per-phase analysis. We describe a per-phase positive-sequence network 
in this subsection and a per-phase analysis procedure in the next subsection.
We make the following simplifying assumptions:
\bi
\AssumptionInc{part:networks; ch:mun.bim; sec:BalancedNk; subsec:1pNetwork; cond:0gammaY}
\item[] C\ref{part:networks; ch:mun.bim; sec:BalancedNk; subsec:1pNetwork; cond:0gammaY}:
    The neutral voltages $\gamma_j := V_j^n = 0$ for all $Y$-configured devices
    $j \in N_v^Y \cup N_c^Y \cup N_i^Y$. 

\AssumptionInc{part:networks; ch:mun.bim; sec:BalancedNk; subsec:1pNetwork; cond:0gammaD}
\item[] C\ref{part:networks; ch:mun.bim; sec:BalancedNk; subsec:1pNetwork; cond:0gammaD}:
    The zero-sequence voltages $\gamma_j := \frac{1}{3}\textbf 1^{\sf T}V_j = 0$ 
    for all $\Delta$-configured voltage sources $j \in N_v^\Delta$. 
\ei
While $\gamma_j$ in 
C\ref{part:networks; ch:mun.bim; sec:BalancedNk; subsec:1pNetwork; cond:0gammaY} and
C\ref{part:networks; ch:mun.bim; sec:BalancedNk; subsec:1pNetwork; cond:0gammaD}
are part of the device specification, the zero-sequence voltages $\gamma_j$ of
$\Delta$-configured current sources and impedances $j \in N_c^\Delta \cup N_i^\Delta$
are not specified but need to be determined through the network equation 
\eqref{part:networks; ch:mun.bim; sec:BalancedNk; subsec:3pAnalysis; eq:I=Y1phiV.1}.
We will prove in Lemma \ref{ch:mun.1; sec:bim; subsec:BalancedNk; lemma:balanced.1}
below that 
C\ref{part:networks; ch:mun.bim; sec:BalancedNk; subsec:1pNetwork; cond:0gammaY} and
C\ref{part:networks; ch:mun.bim; sec:BalancedNk; subsec:1pNetwork; cond:0gammaD}
indeed imply $\gamma_j = 0$ for $j \in N_c^\Delta \cup N_i^\Delta$.
We will explain in Section 
\ref{part:networks; ch:mun.bim; sec:BalancedNk; subsec:extension} how the results
extend to the general case where these assumptions do not hold.

\noindent\textbf{Network equation.}
Consider a network whose graph is
$G = (\overline N, E)$ as before but each line $(j,k)\in E$ is characterized by the complex \emph{scalar} 
admittances $\left( \eta_{jk}^s, \eta_{jk}^m, \eta_{kj}^m \right) \in\mathbb C^3$ in
\eqref{ch:mun.1; sec:bim; subsec:BalancedNk; eq:1a}, instead of $3\times 3$ admittance
matrices in the three-phase network.
Associated with each bus $j\in\overline N$ is a scalar voltage and a scalar current injection
$\left( v_j, i_j \right) \in \mathbb C^2$.  
The current vector $i := (i_j, j\in\overline N)$ and the voltage vector $v := (v_j, j\in\overline N)$ are
related by the $(N+1)\times (N+1)$ per-phase admittance matrix $Y^{1\phi}$ defined in \eqref{eq:defYjk1phi}
according to $i \, = \, Y^{1\phi}v$.
This relationship defines a \emph{per-phase positive-sequence network}.
On this per-phase network, the single-phase devices on buses 
$j\in N_v^\Delta \cup N_c^\Delta \cup N_i^\Delta$ are specified as 
(from Table \ref{ch:mun.1; sec:bim; subsec:BalancedNk; table:buses}):
\begin{align*}
\text{Voltage source: } & & 
v_j & \ = \ \hat\alpha_j \lambda_j, & 	j & \ \in \ N_v
\\
\text{Current source: } &  &
i_j &  \ = \ - \hat\alpha_j \mu_j, 	& 	 j & \ \in \ N_c
\\
\text{Impedance: }  & &
 i_j & \ = \ - \hat\alpha_j \epsilon \, v_j, 	& 	j & \ \in\ N_i
\end{align*}
Note that the per-phase impedance model does not involve $\gamma_j$.
To express this specification in vector form, define the following (column) vectors:
\begin{align*}
\left( v_v , i_v \right) & \ \ := \ \ \left( v_j, i_j, \,  j\in N_v \right)	\\
\left( v_c, i_c \right)& \ \ := \ \ \left( v_j, i_j, \, j\in N_c \right)		\\
\left( v_i, i_i \right) & \ \ := \ \ \left( v_j, i_j, \, j\in N_i \right)
\end{align*}
The voltage sources and impedance modeled are then specified as:
\begin{subequations}
	\begin{align}
	v_v  & \ = \ \left( \hat\alpha_j \lambda_j, \, j\in N_v \right),		\qquad
	i_c   \ = \ \left( - \hat\alpha_j \mu_j, \, j \in N_c \right)
	\\
	i_i & \ = \ - Y_i^{1\phi} v_i  	\ \  \text{ with } \ \
	Y_i^{1\phi} := \diag\left(\hat\alpha_j \epsilon_j, \, j\in N_i \right)
	\label{ch:mun.1; sec:bim; subsec:BalancedNk; eq:1pspec.c}
	\end{align}
	This
	\label{ch:mun.1; sec:bim; subsec:BalancedNk; eq:1pspec}
	\end{subequations}
is the per-phase version of the specification \eqref{ch:mun.1; sec:bim; subsec:BalancedNk; eq:3pspec.1}
of the corresponding three-phase devices.

A key step in per-phase analysis is to solve a per-phase version of the three-phase
problem \eqref{part:networks; ch:mun.bim; sec:BalancedNk; subsec:3pAnalysis; eq:I=Y1phiV.1}\eqref{ch:mun.1; sec:bim; subsec:BalancedNk; eq:3pspec.1}:
given the specification in \eqref{ch:mun.1; sec:bim; subsec:BalancedNk; eq:1pspec}, compute 
the remaining variables
\begin{align*}
v_{-v} & \ \ := \ \ \left( v_j, \, j\not\in N_v \right), 	&    i_{-c}  & \ \ := \ \ \left( i_j, \, j\not\in N_c \right)     & &
\end{align*}
from $i = Y^{1\phi} v$, or equivalently, from
\begin{align}
\begin{bmatrix} i_v \\ i_c  \\  i_i  \end{bmatrix} & \ \ = \ \ 
\underbrace{ \begin{bmatrix} Y^{1\phi}_{vv} & Y^{1\phi}_{vc} & Y^{1\phi}_{vi}  \\  
	Y^{1\phi}_{cv} & Y^{1\phi}_{cc} & Y^{1\phi}_{ci}  \\  
	Y^{1\phi}_{iv} & Y^{1\phi}_{ic} & Y^{1\phi}_{ii}  	\end{bmatrix} }_{Y^{1\phi}}
	\begin{bmatrix} v_v \\ v_c  \\  v_i  \end{bmatrix}  & &
\label{ch:mun.1; sec:bim; subsec:BalancedNk; eq:i=Y1phiv.2}
\end{align}
where the admittance matrix $Y^{1\phi}$ is defined in  \eqref{eq:defYjk1phi}. 
We will first compute $v_{-v}$ from \eqref{ch:mun.1; sec:bim; subsec:BalancedNk; eq:i=Y1phiv.2} and then
substitute $v_{-v}$ back into \eqref{ch:mun.1; sec:bim; subsec:BalancedNk; eq:i=Y1phiv.2} to compute
$i_{-c}$.

\noindent\textbf{Computation of $\left( v_{-v}, i_{-c} \right)$.}
To compute $v_{-v}$ define the following matrices from the per-phase admittance matrix $Y^{1\phi}$:
\begin{subequations}
\begin{align}
\underbrace{ \left[
\begin{array} {c | c c}
	Y^{1\phi}_{vv} & Y^{1\phi}_{vc} & Y^{1\phi}_{vi}  \\  \hline
	Y^{1\phi}_{cv} & Y^{1\phi}_{cc} & Y^{1\phi}_{ci}  \\  
	Y^{1\phi}_{iv} & Y^{1\phi}_{ic} & Y^{1\phi}_{ii}  	\end{array}
\right] }_{Y^{1\phi}}
& \ \ =: \ \ 
\left[ \begin{array}  {c | c} A_{11} & A_{21}^{\sf T}  \\ \hline  A_{21}  & A_{22}  \end{array} \right]	\\
\begin{bmatrix} Y^{1\phi}_{cc} & Y^{1\phi}_{ci} \\ Y^{1\phi}_{ic}  & Y^{1\phi}_{ii} + Y_i^{1\phi} \end{bmatrix}
& \ \ =: \ \ A_{22}'
\label{ch:mun.1; sec:bim; subsec:BalancedNk; thm:1phiNk; eq:Ysubmatrices.1}
\end{align}
where the matrix $Y_i^{1\phi}$ is defined in \eqref{ch:mun.1; sec:bim; subsec:BalancedNk; eq:1pspec.c}.
Note that both $A_{22}$ are $A_{22}'$ are complex symmetric and therefore legitimate admittance 
matrices (they will be interpreted below as admittance matrices of a reduced network
in \eqref{ch:mun.1; sec:bim; subsec:BalancedNk; eq:i=Y1phiv.3}).
For the computation of $i_{-c}$ in 
\eqref{ch:mun.bim; sec:BalancedNk; subsec:1phiAnalysis; eq:i-c}
below, define the following submatrices of $Y^{1\phi}$:
\begin{align}
Y_{-c}^{1\phi} & \ \ := \ \ 
\left[
\begin{array} {c | c c}
	Y^{1\phi}_{vv} & Y^{1\phi}_{vc} & Y^{1\phi}_{vi}  \\ 
	Y^{1\phi}_{iv} & Y^{1\phi}_{ic} & Y^{1\phi}_{ii}  	\end{array}
\right] 
\ \ =: \ \ \left[ \begin{array} {c | c} Y^{1\phi}_{-c,v} & Y^{1\phi}_{-c,-v} \end{array} \right]
\label{ch:mun.1; sec:bim; subsec:BalancedNk; thm:1phiNk; eq:Ysubmatrices.2}
\end{align}
\end{subequations}

Substituting the specification \eqref{ch:mun.1; sec:bim; subsec:BalancedNk; eq:1pspec}
into \eqref{ch:mun.1; sec:bim; subsec:BalancedNk; eq:i=Y1phiv.2} yields a
system of $(|N_c|+|N_i|)$ linear equations in $(|N_c|+|N_i|)$ unknown voltages $v_{-v}$:
	\begin{align*}
	\underbrace{
	\begin{bmatrix} Y^{1\phi}_{cc} & Y^{1\phi}_{ci} \\ Y^{1\phi}_{ic}  & Y^{1\phi}_{ii} + Y_i^{1\phi} \end{bmatrix} }_{A_{22}'}
	\begin{bmatrix} v_c \\ v_i \end{bmatrix}
	& \ \ = \ \ 	
	\underbrace{ \begin{bmatrix} i_c \\ 0 \end{bmatrix} }_{i_{-v}} \ - \
	\underbrace{ \begin{bmatrix} Y^{1\phi}_{cv} \\  Y^{1\phi}_{iv}  \end{bmatrix} }_{A_{21}} v_v
	\end{align*}
or in abbreviation:
\begin{subequations}
\begin{align}
A_{22}' v_{-v} & 	\ \ = \ \   i_{-v} \ - \ A_{21} v_v  \ \ =: \ \ b
\label{ch:mun.1; sec:bim; subsec:BalancedNk; eq:i=Y1phiv.3}
\end{align}
where $i_{-v} := (i_c, 0)$ and the matrices $A_{21}, A_{22}'$ are defined in 
\eqref{ch:mun.1; sec:bim; subsec:BalancedNk; thm:1phiNk; eq:Ysubmatrices.1}.
If $A_{22}'$ is invertible, this yields a unique $v_{-v}$ and hence a unique $i_{-v}$.
Since $A_{22}'$ can be treated as an admittance matrix, the equation 
\eqref{ch:mun.1; sec:bim; subsec:BalancedNk; eq:i=Y1phiv.3} can be interpreted as a $V$-$I$
relationship on a reduced network consisting only of buses in $N_c \cup N_i$ with current sources 
and impedances.  The current sources inject currents $i_c$.  The effect of voltage sources $j\in N_v$ 
on this reduced network is summarized as additional current injections $-A_{21} v_v$ at these buses,
so that the net current injections are $b := i_{-v} - A_{21} v_v$.

Substituting $v_{-v}$ from \eqref{ch:mun.1; sec:bim; subsec:BalancedNk; eq:i=Y1phiv.3}
 into \eqref{ch:mun.1; sec:bim; subsec:BalancedNk; eq:i=Y1phiv.2} and
using \eqref{ch:mun.1; sec:bim; subsec:BalancedNk; thm:1phiNk; eq:Ysubmatrices.2} we have
\begin{align}
i_{-c} & \ \ = \ \ Y_{-c}^{1\phi} v \ \ = \ \ Y_{-c, v}^{1\phi} \, v_v \ + \ Y_{-c, -v}^{1\phi} \left( A_{22}'^{-1} b \right)
\label{ch:mun.bim; sec:BalancedNk; subsec:1phiAnalysis; eq:i-c}
\end{align}
assuming $A_{22}'$ is invertible.
\label{ch:mun.bim; sec:BalancedNk; subsec:1phiAnalysis; eq:vi}
\end{subequations}

\subsection{Per-phase analysis}
\label{part:networks; ch:mun.bim; sec:BalancedNk; subsec:1pAnalysis}

Suppose the matrices $A_{22}$ and $A_{22}'$ are invertible.
The following per-phase analysis procedure offers a simpler alternative to solving the three-phase problem 
\eqref{part:networks; ch:mun.bim; sec:BalancedNk; subsec:3pAnalysis; eq:I=Y1phiV.1}\eqref{ch:mun.1; sec:bim; subsec:BalancedNk; eq:3pspec.1} 
directly:
\bee
\item Solve the positive-sequence network problem 
	\eqref{ch:mun.bim; sec:BalancedNk; subsec:1phiAnalysis; eq:vi} to obtain 
	$v_{-v} := \left( v_j, j\not\in N_v\right)$ and $i_{-c} := \left( i_j, j\not\in N_c \right)$.
	The terminal variables of the three-phase problem are then 
	\begin{align*}
	V & \ \ = \ \ v\otimes \alpha_+,    &  
	I & \ \ = \ \ i\otimes \alpha_+ & & 
	\end{align*}

\item From the terminal voltage and current $(V, I)$,
	the internal voltages, currents and powers of $Y$-configured devices are obtained from the
	conversion rule \eqref{part: NetworkModels; ch:UnbalancedNk; sec:overview; eq:ExtBahevior.1}:
	\begin{align*}
	V_j^Y & \ \ = \ \ v_j\, \alpha_+,   &
	I_j^Y & \ = \ - i_j \, \alpha_+
	\\
	s_j^Y & \ \ = \ \ - v_j \overline{i}_j \textbf 1, 	&	
	j & \ \in \  N_v^Y \cup N_c^Y \cup N_i^Y
	\end{align*}

\item Those of $\Delta$-configured devices can be computed by applying the conversion rule 
	\eqref{ch:mun.1; sec:ComponentModels; subsec:GenLoadD; eq:l2landterminal.s1} to $(V, I)$:
	\begin{align*}
	V_j^\Delta & \ \ = \ \ (1-\alpha) v_j\, \alpha_+, 	&
	I_j^\Delta & \ = \ - \frac{1}{3}\left( 1 - \alpha \right) i_j\, \alpha_+
	\\
	s_j^\Delta & \ \ = \ \ - v_j \overline{i}_j \textbf 1, 	&
		j & \ \in \  N_v^\Delta \cup N_c^\Delta \cup N_i^\Delta
	\end{align*}
\eee
where we have used $\diag\left( \alpha_+ \alpha_+^{\sf H}\right) = \textbf 1$ 
and $\left(1- \alpha \right)\left(1 - \overline{\alpha} \right) = 3$ in computing the internal
variables.

To prove this per-phase analysis procedure solves the three-phase problem 
\eqref{part:networks; ch:mun.bim; sec:BalancedNk; subsec:3pAnalysis; eq:I=Y1phiV.1}\eqref{ch:mun.1; sec:bim; subsec:BalancedNk; eq:3pspec.1} 
when the network is balanced, it suffices to justify Step 1 of the procedure.  This is 
stated in the following theorem.

\begin{theorem}[Per-phase analysis]
\label{ch:mun.1; sec:bim; subsec:BalancedNk; thm:1phiNk.C1C2}
Suppose $A_{22}$, $A_{22}'$ are invertible and 
assumptions 
C\ref{part:networks; ch:mun.bim; sec:BalancedNk; subsec:1pNetwork; cond:0gammaY} and
C\ref{part:networks; ch:mun.bim; sec:BalancedNk; subsec:1pNetwork; cond:0gammaD}
hold.
Let $\left( v, i \right)$ be the unique solution of the positive-sequence network
\eqref{ch:mun.bim; sec:BalancedNk; subsec:1phiAnalysis; eq:vi}.
Then $\left( V, I \right)$ is the unique three-phase solution to
\eqref{part:networks; ch:mun.bim; sec:BalancedNk; subsec:3pAnalysis; eq:I=Y1phiV.1}\eqref{ch:mun.1; sec:bim; subsec:BalancedNk; eq:3pspec.1}
if and only if $V = v \otimes \alpha_+$ and $I = i \otimes \alpha_+$.
\end{theorem}

\subsection{Extension}
\label{part:networks; ch:mun.bim; sec:BalancedNk; subsec:extension}

Theorem \ref{ch:mun.1; sec:bim; subsec:BalancedNk; thm:1phiNk.C1C2} says that
when the voltage and current sources $\left( V_v, I_c \right)$
are balanced positive-sequence sets, all terminal voltages and currents
$(V, I)$ are balanced positive-sequence sets, as long as assumptions 
C\ref{part:networks; ch:mun.bim; sec:BalancedNk; subsec:1pNetwork; cond:0gammaY} and
C\ref{part:networks; ch:mun.bim; sec:BalancedNk; subsec:1pNetwork; cond:0gammaD}
hold.
Without assumptions 
C\ref{part:networks; ch:mun.bim; sec:BalancedNk; subsec:1pNetwork; cond:0gammaY} and
C\ref{part:networks; ch:mun.bim; sec:BalancedNk; subsec:1pNetwork; cond:0gammaD},
Theorem \ref{ch:mun.1; sec:bim; subsec:BalancedNk; thm:1phiNk.C1C2} needs to be modified
to the following.
\begin{theorem}[Per-phase analysis]
\label{ch:mun.1; sec:bim; subsec:BalancedNk; thm:1phiNk}
Suppose $A_{22}$, $A_{22}'$ are invertible.  
Let $\left( v, i \right)$ be the unique solution of the positive-sequence network
\eqref{ch:mun.bim; sec:BalancedNk; subsec:1phiAnalysis; eq:vi}.
Then $\left( V, I \right)$ is the unique three-phase solution to
\eqref{part:networks; ch:mun.bim; sec:BalancedNk; subsec:3pAnalysis; eq:I=Y1phiV.1}\eqref{ch:mun.1; sec:bim; subsec:BalancedNk; eq:3pspec.1}
if and only if
\begin{align*}
V & \ \ = \ \ v \otimes \alpha_+ + \tilde \gamma \otimes \textbf 1
\\
I & \ \ = \ \ i \otimes \alpha_+ \ + \ \tilde \beta \otimes \textbf 1
\end{align*}
for some $\left( \tilde\gamma, \tilde\beta \right) \in \mathbb C^{2(N+1)}$.
\end{theorem}
The theorem says that, without assumptions 
C\ref{part:networks; ch:mun.bim; sec:BalancedNk; subsec:1pNetwork; cond:0gammaY} and
C\ref{part:networks; ch:mun.bim; sec:BalancedNk; subsec:1pNetwork; cond:0gammaD}, 
$(V, I)$ also has a zero-sequence component $\left( \tilde\gamma, \tilde\beta \right)$
in addition to the positive-sequence component $(v, i)$.
While $(v, i)$ is still computed from the per-phase positive-sequence network as described
in Section \ref{part:networks; ch:mun.bim; sec:BalancedNk; subsec:1pNetwork},
$\left( \tilde\gamma, \tilde\beta \right)$ needs to be computed from a
separate per-phase zero-sequence network that is driven by the zero-sequence
voltages $\gamma_v$ of both $Y$ and $\Delta$-configured voltage sources.
The computation of $\left( \tilde\gamma, \tilde\beta \right)$ is more complicated;
see \cite{Low2022}.

Given the terminal voltage and current $(V, I)$, Steps 2 and 3 of the per-phase
analysis procedure in 
Section \ref{part:networks; ch:mun.bim; sec:BalancedNk; subsec:1pAnalysis}
are modified to:
\bee
\item[2. ] From the terminal voltage and current $(V, I)$,
	the internal voltages, currents and powers of $Y$-configured devices are obtained from the
	conversion rule \eqref{part: NetworkModels; ch:UnbalancedNk; sec:overview; eq:ExtBahevior.1} 
	$\left(\gamma_j = V_j^n\right)$: for $j \in N_v^Y \cup N_c^Y \cup N_i^Y$
	\begin{align*}
	V_j^Y  \ \ = \ \ & V_j \ - \ \gamma_j \textbf 1 \ \ = \ \ v_j\, \alpha_+ \ - \ 
		\left( \gamma_j - \tilde\gamma_j \right) \textbf 1
	\\
	I_j^Y  \ \ = \ \ & - I_j \ \ = \ \ - i_j \, \alpha_+ \ + \ \tilde\beta_j \textbf 1
	\\
	s_j^Y 	
		\ \ = \ \ & -\left( v_j i_j^{\sf H} + \left( \gamma_j - \tilde\gamma_j \right) \tilde\beta_j^{\sf H} \right) \textbf 1 
		\\ & 
			+  v_j \tilde\beta_j^{\sf H}\, \alpha_+  \ + \ \left( \gamma_j - \tilde\gamma_j \right) i_j^{\sf H}\, \alpha_-
	\end{align*}

\item[3. ] Those of $\Delta$-configured devices can be computed by applying the conversion rule 
\eqref{ch:mun.1; sec:ComponentModels; subsec:GenLoadD; eq:l2landterminal.s1} to $(V, I)$:
	for $j \in N_v^\Delta \cup N_c^\Delta \cup N_i^\Delta$
	\begin{align*}
	V_j^\Delta & \ \ = \ \ \Gamma V_j \ \ = \ \ (1-\alpha) v_j\, \alpha_+
	\\
	I_j^\Delta & \ \ = \ \ - \Gamma^{\sf T\dag} I_j \ + \ \beta_j \textbf 1 
		\ \ = \ \ - \frac{1}{3}\left( 1 - \alpha \right) i_j\, \alpha_+ \, + \ \beta_j \textbf 1
	\\
	s_j^\Delta & \ \ = \ \  -  v_j i_j^{\sf H} \textbf 1 \ + \ \left( 1 - \alpha\right) v_j \beta_j^{\sf H}\, \alpha_+
	\end{align*}
\eee
In particular, whereas the internal powers $s^Y = v_j i_j^{\sf H}\textbf 1$ 
are identical across the three single-phase devices if assumptions 
C\ref{part:networks; ch:mun.bim; sec:BalancedNk; subsec:1pNetwork; cond:0gammaY} and
C\ref{part:networks; ch:mun.bim; sec:BalancedNk; subsec:1pNetwork; cond:0gammaD}
hold, they are generally different otherwise due to the presence of both
positive and negative-sequence components.

\subsection{Proof of Theorem \ref{ch:mun.1; sec:bim; subsec:BalancedNk; thm:1phiNk.C1C2}}
\label{part:networks; ch:mun.bim; sec:BalancedNk; subsec:equivalence}

Consider $V_{-v} := \left( V_j, j\not\in N_v \right)$.  Recall that 
$\gamma_j := V_j^n$ denote the neutral voltages for $Y$-configured devices 
$j\in N_c^Y \cup N_i^Y$ and $\gamma_j := \frac{1}{3}\textbf 1^{\sf T} V_j$
denote the zero-sequence voltages of $\Delta$-configured devices
$j\in N_c^\Delta \cup N_i^\Delta$.  Under assumption
C\ref{part:networks; ch:mun.bim; sec:BalancedNk; subsec:1pNetwork; cond:0gammaY},
$\gamma_j = 0$ for $j\in N_c^Y \cup N_i^Y$.
The main ingredient of the proof of Theorem \ref{ch:mun.1; sec:bim; subsec:BalancedNk; thm:1phiNk.C1C2}
is the following lemma that implies  that $\gamma_j = 0$ for 
$j\in N_c^\Delta \cup N_i^\Delta$ as well and that $V_{-v}$ is balanced.
\begin{lemma}[Balanced voltage]
\label{ch:mun.1; sec:bim; subsec:BalancedNk; lemma:balanced.1}
Suppose $(V, I)$ satisfies \eqref{part:networks; ch:mun.bim; sec:BalancedNk; subsec:3pAnalysis; eq:I=Y1phiV.1}\eqref{ch:mun.1; sec:bim; subsec:BalancedNk; eq:3pspec.1}.
Under assumptions 
C\ref{part:networks; ch:mun.bim; sec:BalancedNk; subsec:1pNetwork; cond:0gammaY} and
C\ref{part:networks; ch:mun.bim; sec:BalancedNk; subsec:1pNetwork; cond:0gammaD}:
\bee
\item If $A_{22}$ is invertible, then $\gamma_j = 0$ for 
    $j\in N_c^\Delta \cup N_i^\Delta$ and hence $\gamma_j = 0$ for all $j\in\overline N$.
    
\item The terminal voltage $V = (V_{v}, V_{-v})$ satisfies
    \begin{align*}
	\left( A_{22}' \otimes \mathbb I \right) V_{-v}  & \ \ = \ \ 
	    \left( i_{-v} - A_{21}v_v \right) \otimes \alpha_+  
	    \ \ =: \ \ b \otimes \alpha_+ 
	\end{align*}
	where 
    $A_{22}'$ is defined in 
	\eqref{ch:mun.1; sec:bim; subsec:BalancedNk; thm:1phiNk; eq:Ysubmatrices.1}.
     If $A_{22}'$ is invertible, then $V_{-v} \, = \, \left( A_{22}'^{-1}b \right) \otimes \alpha_+$ is balanced.
\eee
\end{lemma}
\begin{proof}[Proof of Lemma \ref{ch:mun.1; sec:bim; subsec:BalancedNk; lemma:balanced.1}]
For part 1 of the lemma, observe from 
Table \ref{ch:mun.1; sec:bim; subsec:BalancedNk; table:buses} that 
$\textbf 1^{\sf T} I_j = 0$ for $j\not\in N_v^Y$. 
Define $\hat \gamma_j \ := \ \frac{1}{3}\textbf 1^{\sf T} V_j$.  Note that
$\hat\gamma_j = \gamma_j$ for $\Delta$-configured devices but not necessarily 
for $Y$-configured current sources or impedances.
Multiplying both sides of  \eqref{ch:mun.1; sec:bim; subsec:BalancedNk; eq:1b} by 
$\textbf 1^{\sf T}$ gives 
\begin{align*}
\sum_{k: j\sim k}\left( y^s_{jk} \ + \ {y_{jk}^m} \right) \hat\gamma_j \ - 
\sum_{k: j\sim k} y^s_{jk} \hat\gamma_k	& \ \ = \ \ 
\left\{\begin{array}{lcl}
	\frac{1}{3} \textbf 1^{\sf T} I_j	& \text{ if } & j\in N_v^Y	\\	0	& \text{ if } & j\not\in N_v^Y 
\end{array} \right.
& &	
\end{align*}
We can write this for $j\not\in N_v$ in vector form using 
\eqref{ch:mun.1; sec:bim; subsec:BalancedNk; thm:1phiNk; eq:Ysubmatrices.1} as:
\begin{align*}
\left[ \begin{array}  {c | c}  A_{21}  & A_{22}  \end{array} \right] 
\begin{bmatrix} \hat\gamma_v \\ \hat\gamma_{-v} \end{bmatrix}
\ \ = \ \ 0
\end{align*}
Note that $\hat\gamma_v = \gamma_v$ because for $j\in N_v^Y$, 
$\hat \gamma_j := \frac{1}{3} \textbf 1^{\sf T} \left( V_j^Y + V_j^n \textbf 1 \right) = V_j^n = \gamma_j$ and for $j\in N_v^\Delta$,
$\hat\gamma_j = \gamma_j$ by definition.
This implies that $A_{21} \gamma_v + A_{22} \hat\gamma_{-v} \ = \ 0$.
When $A_{22}$ is invertible,
\begin{align*}
\hat\gamma_{-v} & \ \ = \ \ - A_{22}^{-1} A_{21} \gamma_v
\end{align*}
Under assumptions 
C\ref{part:networks; ch:mun.bim; sec:BalancedNk; subsec:1pNetwork; cond:0gammaY} and
C\ref{part:networks; ch:mun.bim; sec:BalancedNk; subsec:1pNetwork; cond:0gammaD},
$\gamma_v = 0$ and hence $\hat\gamma_{-v} = 0$.
This implies $\gamma_j = \hat\gamma_j = 0$ for $j\in N_c^\Delta \cup N_i^\Delta$
and hence $\gamma_j = 0$ for all $j\in\overline N$.
This establishes part 1 of the lemma.

We will prove part 2 in three steps.  First we will express the three-phase specification 
\eqref{ch:mun.1; sec:bim; subsec:BalancedNk; eq:3pspec.1} in terms of per-phase variables.
Then we will write the three-phase equation  
\eqref{part:networks; ch:mun.bim; sec:BalancedNk; subsec:3pAnalysis; eq:I=Y1phiV.1} 
using per-phase parameters.
Finally we will establish the equation for $V_{-v}$ in part 2 of the lemma. 

First, using \eqref{ch:mun.1; sec:bim; subsec:BalancedNk; eq:1pspec}, 
the three-phase quantities in
\eqref{ch:mun.1; sec:bim; subsec:BalancedNk; eq:3pspec.1} can be written as
(under C\ref{part:networks; ch:mun.bim; sec:BalancedNk; subsec:1pNetwork; cond:0gammaY} and
C\ref{part:networks; ch:mun.bim; sec:BalancedNk; subsec:1pNetwork; cond:0gammaD}):
\begin{subequations}
\begin{align}
V_v & \ \ = \ \ \left( \hat\alpha_j \lambda_j \alpha_+ + \gamma_j\textbf 1, \ j \in N_v \right)
	\ \ = \ \ v_v \otimes \alpha_+ 
\\
I_c & \ \ = \ \ \left( - \hat\alpha_j \mu_j \alpha_+, \ j \in N_c \right) \ \ = \ \ i_c \otimes \alpha_+
\\
I_i & \ \ = \ \ \left( - \hat\alpha_j \epsilon_j \left( V_j - \gamma_j \textbf 1 \right), \ j \in N_i \right)
	\ \ = \ \ - \left( Y_i^{1\phi} \otimes \mathbb I \right) V_i 
\end{align}
where
\label{ch:mun.1; sec:bim; subsec:BalancedNk; eq:3pspec1p}
\end{subequations}
the matrix $Y_i^{1\phi}$ is defined in \eqref{ch:mun.1; sec:bim; subsec:BalancedNk; eq:1pspec.c}.

Second, write \eqref{part:networks; ch:mun.bim; sec:BalancedNk; subsec:3pAnalysis; eq:I=Y1phiV.1} in terms
of submatrices of $Y^{1\phi}$ in \eqref{ch:mun.1; sec:bim; subsec:BalancedNk; eq:i=Y1phiv.2}
(since $Y = Y^{1\phi}\otimes \mathbb I$):
\begin{align*}
\begin{bmatrix} I_v \\ I_c  \\  I_i  \end{bmatrix} & \ \ = \ \ 
\underbrace{ \begin{bmatrix} Y_{vv} & Y_{vc} & Y_{vi}  \\  Y_{cv} & Y_{cc} & Y_{ci}  \\  Y_{iv} & Y_{ic} & Y_{ii}  
	\end{bmatrix} }_{Y}
\begin{bmatrix} V_v \\ V_c  \\  V_i  \end{bmatrix}
\\
& \ \ = \ \ 
\left(
\underbrace{ \begin{bmatrix} Y^{1\phi}_{vv} & Y^{1\phi}_{vc} & Y^{1\phi}_{vi}  \\  Y^{1\phi}_{cv} & Y^{1\phi}_{cc} & Y^{1\phi}_{ci}  \\  Y^{1\phi}_{iv} & Y^{1\phi}_{ic} & Y^{1\phi}_{ii}  
	\end{bmatrix} }_{Y^{1\phi}} \otimes\ \mathbb I
\right)	\begin{bmatrix} V_v \\ V_c  \\  V_i  \end{bmatrix}
\end{align*}
Therefore the voltages $(V_c, V_i)$ satisfies
\begin{align*}
I_c & \ \ = \ \ \left( Y^{1\phi}_{cv} \otimes \mathbb I \right) V_v \ + \ 
	\left( Y^{1\phi}_{cc} \otimes \mathbb I \right) V_c \ + \ 
	\left( Y^{1\phi}_{ci} \otimes \mathbb I \right) V_i
\\
I_i & \ \ = \ \ \left( Y^{1\phi}_{iv} \otimes \mathbb I \right) V_v \ + \ 
	\left( Y^{1\phi}_{ic} \otimes \mathbb I \right) V_c \ + \ 
	\left( Y^{1\phi}_{ii} \otimes \mathbb I \right) V_i
\end{align*}
Substitute \eqref{ch:mun.1; sec:bim; subsec:BalancedNk; eq:3pspec1p} to get
\begin{align*}
\left( Y^{1\phi}_{cc} \otimes \mathbb I \right) V_c +
\left( Y^{1\phi}_{ci} \otimes \mathbb I \right) V_i		 \ = \ &
i_c \otimes \alpha_+ -
\left( Y^{1\phi}_{cv} \otimes \mathbb I \right) \left( v_v \otimes \alpha_+ \right)
\\
\left( Y^{1\phi}_{ic} \otimes \mathbb I \right) V_c + 
\left( Y^{1\phi}_{ii} \otimes \mathbb I \right) V_i		 \  = \ &
- \left( Y_i^{1\phi} \otimes \mathbb I \right) V_i
\\	& 
- \left( Y^{1\phi}_{iv} \otimes \mathbb I \right) \left( v_v \otimes \alpha_+ \right)
\end{align*}
Since 
$\left( Y^{1\phi}_{cv} \otimes \mathbb I \right) \left( v_v \otimes \alpha_+ \right) = 
\left( Y^{1\phi}_{cv} v_v \right) \otimes \alpha_+$ and 
$\left( Y^{1\phi}_{iv} \otimes \mathbb I \right) \left( v_v \otimes \alpha_+ \right) = 
\left( Y^{1\phi}_{iv} v_v \right) \otimes \alpha_+$, the voltage $(V_c, V_i)$ satisfies
\begin{align*}
\left(
\underbrace{ 
\begin{bmatrix}	Y^{1\phi}_{cc} & Y^{1\phi}_{ci} \\  Y^{1\phi}_{ic}  &   Y^{1\phi}_{ii} + Y_i^{1\phi}
\end{bmatrix} }_{A_{22}'}
	\otimes \ \mathbb I	\right)
\begin{bmatrix} V_c \\ V_i  \end{bmatrix}	& \ = \ \ 
\underbrace{
\left(
\begin{bmatrix}  i_c  \\ 0 \end{bmatrix} - A_{21} v_v \right) }_{b} \otimes \ \alpha_+
\end{align*}
This is the equation in part 2 of the lemma, abbreviated as
\begin{align}
\left( A_{22}' \otimes \mathbb I \right) V_{-v}  & \ \ = \ \ b \otimes \alpha_+ 
\label{ch:mun.1; sec:bim; subsec:BalancedNk; eq:VcViKronecker}
\end{align}
If $A_{22}'$ is invertible, so is $A_{22}' \otimes \mathbb I$.
Then $V_{-v}$ is uniquely determined as
\begin{align*}
V_{-v} & \ \ = \ \ \left( A_{22}'\otimes \mathbb I \right)^{-1} \left( b \otimes \alpha_+  \right)
	 \ \ = \ \ \left( A_{22}'^{-1} \otimes \mathbb I \right) \left( b \otimes \alpha_+ \right)
\\ &
 \ \ = \ \ \left( A_{22}'^{-1}b \right) \otimes \alpha_+ 
\end{align*} 
This completes part 2 of the lemma.
\end{proof}
\vspace{0.2in}

\begin{proof}[Proof of Theorem \ref{ch:mun.1; sec:bim; subsec:BalancedNk; thm:1phiNk.C1C2}]
Under the assumption that $A_{22}$ and $A_{22}'$ are invertible, the solution
$v_{-v}$ to the per-phase network
\eqref{ch:mun.1; sec:bim; subsec:BalancedNk; eq:i=Y1phiv.3} is unique.
Lemma \ref{ch:mun.1; sec:bim; subsec:BalancedNk; lemma:balanced.1} implies that $V_{-v}$ is the
unique solution to \eqref{ch:mun.1; sec:bim; subsec:BalancedNk; eq:VcViKronecker}.
We now show that the expression $V_{-v} = v_{-v} \otimes \alpha_+$ in the 
theorem satisfies this equation and therefore must be its unique solution.  
We have
\begin{align*}
\left( A_{22}' \otimes \mathbb I \right) \left( v_{-v} \otimes \alpha_+ \right)
 & \ \ = \ \ 
\left( A_{22}' v_{-v} \right) \otimes \alpha_+ 
\end{align*}
But $A_{22}' v_{-v} =  i_{-v} - A_{21}v_{v} = b$ from \eqref{ch:mun.1; sec:bim; subsec:BalancedNk; eq:i=Y1phiv.3}.
This proves the expression for $V_{-v}$ in the theorem.

To prove the expression for $I_{-c}$ in the theorem, 
let $V_{-c} := \left( V_v, V_i \right)$.  Substituting into 
\eqref{part:networks; ch:mun.bim; sec:BalancedNk; subsec:3pAnalysis; eq:I=Y1phiV.1} we have
(since $Y = Y^{1\phi}\otimes \mathbb I$)
\begin{align*}
I_{-c} & \ \ = \ \ 
\left(
\underbrace{ \begin{bmatrix} Y^{1\phi}_{vv} & Y^{1\phi}_{vc} & Y^{1\phi}_{vi}  \\ Y^{1\phi}_{iv} & Y^{1\phi}_{ic} & Y^{1\phi}_{ii}  
	\end{bmatrix} }_{Y_{-c}^{1\phi}} \otimes\ \mathbb I
\right)	\begin{bmatrix} V_v \\ V_c  \\  V_i  \end{bmatrix}
	\ \ =: \ \ Y_{-c}^{1\phi} \, V
\end{align*}
where $Y_{-c}^{1\phi}$ is the matrix defined in 
\eqref{ch:mun.1; sec:bim; subsec:BalancedNk; thm:1phiNk; eq:Ysubmatrices.2}.
Substitute the expression $V = v\otimes \alpha_+$ from 
Theorem \ref{ch:mun.1; sec:bim; subsec:BalancedNk; thm:1phiNk.C1C2} we have
\begin{align*}
I_{-c} & 
\ \ = \ \ \left( Y_{-c, v}^{1\phi} v_v + Y_{-c, -v}^{1\phi} v_{-v} \right) \otimes \alpha_+ 
\ \ = \ \ i_{-c} \otimes \alpha_+
\end{align*}
where the last equality follows from $v_{-v} = A_{22}'^{-1}b$ and 
\eqref{ch:mun.bim; sec:BalancedNk; subsec:1phiAnalysis; eq:i-c}.
This completes the proof
of Theorem \ref{ch:mun.1; sec:bim; subsec:BalancedNk; thm:1phiNk.C1C2}.
\end{proof} \vspace{0.2in}

\section{Conclusion}
\label{sec:conc}

In this paper we have introduced a three-phase power flow model.  Its key feature is to model
a three-phase device by (i) an internal model that describes the behavior of the constituent
single-phase devices, and (ii) a conversion rule that maps internal variables to terminal 
variables observable externally of the three-phase device based on its configuration.  
The flexibility allows a unified
approach to compose general network models where devices can be ungrounded, grounded
with zero or nonzero earthing impedances, or have nonzero zero-sequence voltages.
We have illustrated our model by studying a general analysis problem when the network is 
balanced, formalizing per-phase analysis, and proving its validity using the spectral properties
of the linear conversion matrix $\Gamma$.

\section*{Acknowledgment}

We would like to thank 
	Janusz Bialek,
	Frederik Geth,
	Pierluigi Mancarella,
and 
	Danny Tsang 
for  helpful discussions.

\bibliographystyle{unsrt}
\bibliography{PowerRef-201202}

\end{document}